\documentclass[12pt, reqno]{amsart}

\usepackage{marginnote}
\usepackage{geometry}

\usepackage{amsthm}
\usepackage{amsfonts}
\usepackage{amsmath}
\usepackage{amssymb}
\usepackage{mathrsfs}
\usepackage{epsfig}

\makeatletter
\@namedef{subjclassname@2020}{%
  \textup{2020} Mathematics Subject Classification}
\makeatother

\newcommand{\overbar}[1]{\mkern 1.5mu\overline{\mkern-1.5mu#1\mkern-1.5mu}\mkern 1.5mu}

\usepackage{fourier}

\usepackage{xcolor}

 \usepackage{float}

\newcommand{\ncom}{\newcommand}
\ncom{\ol}{\overline}
\ncom{\nno}{\nonumber}
\ncom{\rar}{\rightarrow}
\ncom{\Rar}{\Rightarrow}
\ncom{\noin}{\noindent}
\ncom{\sz}{\scriptsize}
\ncom{\rf}{\ref}
\ncom{\sgm}{\sigma}
\ncom{\Sgm}{\Sigma}
\ncom{\dt}{\delta}
\ncom{\Dt}{Delta}

\ncom{\s}{\underline{s}}
\ncom{\lmd}{\lambda}
\ncom{\Lmd}{\Lambda}
\ncom{\eps}{\epsilon}
\ncom{\pcc}{\stackrel{P}{>}}
\ncom{\dist}{{\rm\,dist}}
\ncom{\sspan}{{\rm\,span}}
\ncom{\re}{{\rm Re\,}}
\ncom{\im}{{\rm Im\,}}
\ncom{\sgn}{{\rm sgn\,}}
\ncom{\eop}{\hfill{{\rule{2.5mm}{2.5mm}}}}
\ncom{\eoe}{\hfill{{\rule{1.5mm}{1.5mm}}}}
\ncom{\eof}{\hfill{{\rule{1.5mm}{1.5mm}}}}
\ncom{\hone}{\mbox{\hspace{1em}}}
\ncom{\htwo}{\mbox{\hspace{2em}}}
\ncom{\hthree}{\mbox{\hspace{3em}}}
\ncom{\hfour}{\mbox{\hspace{4em}}}
\ncom{\hsev}{\mbox{\hspace{7em}}}
\ncom{\vone}{\vskip 2ex}
\ncom{\cH}{{\mathcal H}}
\ncom{\vtwo}{\vskip 4ex}
\ncom{\vonee}{\vskip 1.5ex}
\ncom{\vthree}{\vskip 6ex}
\ncom{\vfour}{\vspace*{8ex}}
\ncom{\norm}{\|\;\;\|}
\ncom{\integ}[4]{\int_{#1}^{#2}\,{#3}\,d{#4}}
\ncom{\inp}[2]{\langle{#1},\,{#2} \rangle}
\ncom{\Inp}[2]{\big\langle{#1},\,{#2} \big\rangle}
\ncom{\vspan}[1]{{{\rm\,span}\#1 \}}}
\ncom{\dm}[1]{\displaystyle {#1}}
\ncom{\Hom}{\operatorname{Hom}}
\ncom{\Hol}{\operatorname{Hol}}
\ncom{\Ps}{\mathcal P_{\underline{s}}}

\ncom{\defin} {\overset {\text {\rm def} }{=}}

\usepackage{soul}
\setstcolor{red}

\newtheorem{thm}{Theorem}[section]

\newtheorem{cor}[thm]{Corollary}
\newtheorem{lem}[thm]{Lemma}

\newtheorem*{lem2.3}{Lemma 2.3}

\newtheorem{prop}[thm]{Proposition}
\theoremstyle{definition}
\newtheorem{defn}[thm]{Definition}
\theoremstyle{remark}
\newtheorem{rem}[thm]{Remark}

\numberwithin{equation}{section}
\newtheorem{conj}[thm]{Conjecture}

\def \s{\underline{s}}
\renewcommand{\epsilon}{\varepsilon}
\renewcommand{\kappa}{\varkappa}

\numberwithin{equation}{section} 

\begin{document}

\title[A trace inequality for commuting tuples of operators]{A trace inequality for commuting tuple of operators}

\author[G. Misra]{Gadadhar Misra}
\address[G. Misra]{Department of Mathematics, Indian Institute of Science, Bangalore 560012, India}
\email[G. Misra]{gm@iisc.ac.in}


\author[P. Pramanick]{Paramita Pramanick}
\address[P. Pramanick]{Department of Mathematics, Indian Institute of Science,
	Bangalore 560012, India} 
	\email[P. Pramanick]{paramitapramanick@gmail.com}

\author[K. B. Sinha]{KALYAN B. Sinha}
\address[K. B. Sinha]{Jawaharlal
Nehru Centre for Advanced Scientific Research (JNCASR), Jakkur,
Bangalore - 560 064, India}
\email[K. B. Sinha]{kbs@jncasr.ac.in}

\thanks{The first-named author would like to acknowledge funding through the J C Bose National Fellowship and the MATRICS grant of SERB. The second named author was supported by a Research Fellowship of the National Board for Higher Mathematics, DAE. The third named author would like to acknowledge the funding he has received through the Senior Scientist program of INSA}
\thanks{A number of the results presented in this paper are from the PhD thesis of the second named author submitted to the Indian Institute of Science, Bangalore}


\date{}

\begin{abstract}
For a commuting $d$- tuple of operators $\boldsymbol T$ defined on a complex separable Hilbert space $\mathcal H$, let $\big [ \!\!\big [ \boldsymbol T^*, \boldsymbol T  \big ]\!\!\big ]$ be the $d\times d$ block operator 
$\big (\!\!\big (\big [  T_j^* , T_i\big ]\big )\!\!\big )$  
of the commutators $[T^*_j , T_i] := T^*_j T_i - T_iT_j^*$. We define the determinant of $\big [ \!\!\big [ \boldsymbol T^*, \boldsymbol T  \big ]\!\!\big ]$ by symmetrizing the products in the Laplace formula for the determinant of a scalar matrix. We prove that the determinant of $\big [ \!\!\big [ \boldsymbol T^*, \boldsymbol T  \big ]\!\!\big ]$ equals the generalized commutator of the $2d$ - tuple of operators, $(T_1,T_1^*, \ldots, T_d,T_d^*)$ introduced earlier by Helton and Howe. We then apply the  Amitsur-Levitzki theorem to conclude that for any commuting $d$ - tuple of $d$ - normal operators, the determinant of $\big [ \!\!\big [ \boldsymbol T^*, \boldsymbol T  \big ]\!\!\big ]$ must be $0$.  We show that if the $d$- tuple $\boldsymbol T$ is cyclic, the determinant of $\big [ \!\!\big [ \boldsymbol T^*, \boldsymbol T  \big ]\!\!\big ]$ is non-negative and the compression of a fixed set of words in $T_j^* $ and $T_i$ -- to a nested sequence of finite dimensional subspaces increasing to $\mathcal H$ -- does not grow very rapidly, then the trace of the determinant of the operator $\big [\!\! \big [ \boldsymbol T^* , \boldsymbol T\big ] \!\!\big  ]$ is finite. Moreover, an upper bound for this trace is given. 
This upper bound is shown to be sharp for a class of commuting $d$ - tuples. We make a conjecture of what might be a sharp bound in much greater generality and verify it in many examples. \end{abstract}

%
%
%
%

\maketitle


\section{Introduction}
The study of bounded linear operators $T$ on a complex separable Hilbert space close to normal operators has been quite successful and has a long history. For commuting tuples of operators $\boldsymbol T:=(T_1, \ldots ,T_d)$, this study is not as complete as in the case of one variable. For instance, while some attempts at the generalization of the Berger-Shaw or the Brown, Douglas and Fillmore theorem in the multi-variate context exist, they are far from complete. In this paper we focus on obtaining trace estimates for the commutators $[T_j^*,T_i]$, $1\leq i,j\leq d$, in the spirit of the Berger-Shaw theorem. As in the one variable case, one expects a reasonable estimate after requiring the $d$ - tuple to be hyponormal and $m$ - cyclic. We recall these notions below. 


An operator $T$ on a Hilbert space $\mathcal{H}$ is said to be \emph{hyponormal} if the commutator $[T^*,T]$ is non-negative definite.  There are many possible notions of hyponormality for a $d$- tuple $\boldsymbol T=(T_1,\ldots,T_d)$ of commuting operators acting on a Hilbert space $\mathcal{H}$. 

For instance, a commuting $d$- tuple $\boldsymbol T$ of operators acting on a Hilbert space $\mathcal{H}$  is said to be  \emph{jointly hyponormal} if $$\big[\!\! \big [\boldsymbol{T}^*, \boldsymbol T\big ]\!\!\big ]:=\big ( \!\! \big ([T_j^*, T_i]\big ) \!\!\big )_{i,j=1}^d:\bigoplus_d\mathcal{H} \longrightarrow \bigoplus_d\mathcal{H}$$
is non-negative definite, that is, for each $x \in \bigoplus_d\mathcal{H}$,
$$ \big \langle \big [\!\!\big [\boldsymbol{T}^*, \boldsymbol T\big ]\!\!\big ]x, x\big \rangle \geq 0.$$


On the other hand, a commuting $d$- tuple $\boldsymbol{T}$  of operators acting on a Hilbert space $\mathcal{H}$ is said to be \emph{projectively hyponormal} if, for each vector $(\alpha_{1},\ldots, \alpha_{d})\in {\mathbb{C}^d}\setminus \{\boldsymbol 0\},$ the sum  $\sum_{i=1}^{d}\alpha_{i}T_{i}$ is a hyponormal operator on $\mathcal{H}.$ In the literature, these have been called weakly hyponormal operators. However, the adjective "projective" describes better the way this class is related to the class of jointly hyponormal operators.  




\begin{defn}[Schatten $p$ - class] \label{defn1.1}
An operator $T:\mathcal{H}\to \mathcal{H}$ is said to be in the Schatten $p$ - class $\mathcal B_p$, $p\in [1,\infty)$, if there is an orthonormal basis $\{e_n\}$ such that
$$\|T\|_p:=\big ( \sum_n  \langle |T|^pe_n, e_n\rangle \big )^{1/p} < \infty,$$
where $|T|$ is the unique positive square root of the operator $T^*T.$
The operator $T$ is said to be of trace class, i.e., $T\in \mathcal B_1$ if $\|T\|_1 < \infty$. We let $\mathcal B(\mathcal H)$ denote the algebra of all bounded operators on the Hilbert space $\mathcal H$. 
\end{defn}

($m$- cyclic). The $d$ - tuple $\boldsymbol T$ is said to be $m$ - polynomially cyclic if the smallest cardinality of a linearly independent set of vectors $\boldsymbol \xi\subseteq \mathcal H$ such that $\mathcal H$ is the closed linear span of
\[\bigg\{T_{1}^{i_1}T_{2}^{i_2}\ldots T_{d}^{i_d}v |~ v\in \boldsymbol{\xi} ~\text{and}~ i_{1}, i_{2}, \ldots ,i_{d}\geq 0  \bigg \}  \]
is equal to $m$. We let $\boldsymbol \xi[m]$ denote any such set $\boldsymbol \xi$ of $m$ linearly independent vectors.

In what follows, we will use the terms ``finitely polynomially cyclic'',  ``$m$ - polynomially cyclic"  or even  "$m$ - cyclic" interchangeably, while ``cyclic'' always means $1$ - cyclic. Also, unless explicitly indicated otherwise,  we  drop  the  adjective  ``jointly''  and  write  ``hyponormal'' instead of jointly hyponormal. 

Let $\sigma(\boldsymbol T)\subseteq \mathbb C^d$ denote the Taylor joint spectrum of the $d$-tuple of operators $\boldsymbol T$, see \cite{Cspectrum, EPut}. We assume that the spectrum  is \textit{polynomially convex}. Consequently, the property of polynomial cyclicity as opposed to rational cyclicity, would be the natural
hypothesis for us. Therefore, we write \textit{``$m$-cyclic''} instead of $m$-polynomially cyclic throughout this paper. 

A second consequence of the polynomial convexity is that $\boldsymbol T$ is $m$- cyclic if and only if the subspace $$ \big\{ f (\boldsymbol T) v |~ v\in \boldsymbol{\xi}[m] ~\text{and}~ f\in \text{Hol}(\sigma(\boldsymbol T) )\big \}$$
is dense in $\mathcal H$. Here, $\text{Hol}(\sigma(\boldsymbol T) )$ is the algebra of all functions which are holomorphic in some open neighbourhood of the compact set $\sigma(\boldsymbol T)$. The operator $f(\boldsymbol T)$ is then defined using the usual holomorphic functional calculus. 

For a single hyponormal operator $T$, one of the surprising results is the trace inequality due to Berger and Shaw reproduced below. A generalization of the Berger-Shaw inequality is in \cite{Voiculescu}, see also \cite{HN}. 

\begin{thm}[Berger-Shaw, \cite{Bergershaw}]
If $T$ is an $m$- cyclic hyponormal operator, then $[T^*, T]$ is in trace-class and
$$\text{trace}\,[T^*, T]\leq \frac{m}{\pi}\nu\,\big (\sigma(T)\big ),$$
where $\sigma(T)$ is the spectrum of $T$ and $\nu$ is the Lebesgue measure. 
\end{thm}
Apart from the earlier results for a $d$ - tuple of $m$ - cyclic hyponormal operators in \cite{Atjoint} and \cite{Douglasyan}, 
more recently, several properties involving commutators of the restriction of the  multiplication by the coordinate functions on $H^2(\mathbb D^2)$
to an invariant subspace and the compression to a co-invariant subspace has been obtained in \cite{Yang}.  

Let $\boldsymbol M:=(M_1,M_2)$ denote the pair of multiplication by coordinate functions $z_1$ and $z_2$, respectively, on the Hardy space $H^2(\mathbb D^2)$ of the bidisc $\mathbb D^2$. The constant vector $1$ is a cyclic vector for $\boldsymbol M$ and  
$$\big[\!\! \big [\boldsymbol{M}^*, \boldsymbol M\big ]\!\!\big ] =  \begin{pmatrix} P\otimes I& 0\\
  0 & I\otimes P \end{pmatrix},$$
where $P$ is the projection to the one dimensional subspace generated by the constant function $1$ in $H^2(\mathbb D)$, is non-negative definite.
Therefore the pair $\boldsymbol M$ is hyponormal. However, neither $[M_1^*,M_1]$ nor $[M_2^*, M_2]$ is trace class. These operators are not even compact. Thus the obvious generalization of the Berger-Shaw inequality fails, in general, in the multivariate case.
However, the product of the two operators $[M_1^*,M_1]$ and $[M_2^*, M_2]$ is trace class with trace norm equal to one. This example prompts the question: In general, for an appropriate class of commuting $d$ - tuple  $\boldsymbol T$, if there is some function of $\big[\!\! \big [\boldsymbol{T}^*, \boldsymbol T\big ]\!\!\big ]$ which might be in the trace class. We identify such a class of commuting $d$ - tuples and a function: $\text{dEt}: \mathcal B(\mathcal H)\otimes \mathcal M_d \to \mathcal B(\mathcal H)$. On the set 
$$\big \{ \big ( \!\!\big ( [T_j^*, T_i ]\big ) \!\!\big )_{i,j=1}^d: T_jT_i = T_iT_j, 1\leq i, j \leq d \big \} \subset \mathcal B(\mathcal H) \otimes \mathcal M_d,$$ it is given by the formula 
$$\text{dEt}\,\big (\big[\!\! \big [\boldsymbol{T}^*, \boldsymbol T\big ]\!\!\big ]\big ):=\sum_{\sigma, \tau \in \mathfrak S_d} \text{Sgn}(\sigma){[T^*_{\tau(1)},T_{\sigma(\tau(1))}]}{[T^*_{\tau(2)},T_{\sigma(\tau(2))}}] \cdots [{T^*_{\tau(d)},T_{\sigma(\tau(d))}}],$$
where $\mathfrak S_d$ is the permutation group on $d$ symbols. 
For the pair of multiplication operators on $H^2(\mathbb D^2)$, $\text{trace}\,\big(\text{dEt}\,\big (\big[\!\! \big [\boldsymbol{M}^*, \boldsymbol M\big ]\!\!\big ]\big )\big) = 2.$

In Section 2, for any hyponormal (resp. weakly hyponormal) $d$ - tuple $\boldsymbol T$, we obtain an estimate, Theorems \ref{Afirst} and \ref{Asecond}, for the trace norm of $[T_j^*,T_i]$, $1\leq i,j\leq d$. Also, we show that  \cite[Theorem 2]{Douglasyan} is valid with the weaker hypothesis of projective hyponormality instead of the hyponrmality assumption. This is Theorem \ref{12345}. In several explicit examples, we verify  the assumptions of this theorem. 

In Section 3, we establish the equality of $\text{dEt}\,\big (\big[\!\! \big [\boldsymbol{T}^*, \boldsymbol T\big ]\!\!\big ]\big )$ with the generalized commutator $\text{GC}(\boldsymbol T^*,\boldsymbol T)$, Definition \ref{Hel-Ho}, introduced earlier by Helton and Howe in Theorem \ref{HH}.  

A commuting $d$ - tuple $\boldsymbol T=(T_1,\ldots , T_d)$  is said to be $n$ - normal if  the operators $T_k$, $1\leq k \leq d$, are from a  $C^*$ - algebra $\mathfrak A$ that is $*$ - isometrically isomorphic to the $C^*$ - algebra $C^*(\hat{\boldsymbol T})\otimes \mathcal M_n$ for some commuting $d$- tuple $\hat{\boldsymbol T}$ of normal operators.

Applying the Amitsur-Levitzki theorem, which states that ``the algebra of $d\times d$ matrices over any commutative ring $R$ satisfies the standard polynomial $S_{2d}$'', see \cite{P}, to a commuting tuple $\boldsymbol T$ of $d$ - normal operators and using the equality of dEt and GC, we show that if  $\boldsymbol T$ is a $d$ - tuple of $d$ - normal operators, then $\text{dEt}\,\big (\big[\!\! \big [{\boldsymbol{T}}^*, {\boldsymbol{T}}\big ]\!\!\big ]\big )=0$. This is Corollary \ref{cor:3.6}. 

In Section 4, We obtain a trace estimate for $\text{dEt}\,\big (\big[\!\! \big [\boldsymbol{T}^*, \boldsymbol T\big ]\!\!\big ]\big )$ assuming that the $d$ - tuple $\boldsymbol T$ is $m$ - cyclic, is in the class $B_{m,\theta}(\Omega)$ that is introduced in this paper and is such that 
$\text{dEt}\,\big (\big[\!\! \big [\boldsymbol{T}^*, \boldsymbol T\big ]\!\!\big ]\big )$ is non-negative definite, see Theorem \ref{lem63}. 

In Section 5, we discuss the weighted Bergman spaces on the Euclidean ball and the ellipsoid. For the commuting tuple $\boldsymbol M$ of multiplication by the coordinate functions $M_i$, $1\leq i \leq d$, in these examples, we compute the trace of the operator  $\text{dEt}\,\big (\big[\!\! \big [\boldsymbol{M}^*, \boldsymbol M\big ]\!\!\big ]\big )$. Based on these computations, we conjecture an inequality for the trace of the operator  $\text{dEt}\,\big (\big[\!\! \big [\boldsymbol{T}^*, \boldsymbol T\big ]\!\!\big ]\big )$ in general, which we expect to be sharp. 


\section{A class of hyponormal $d$- tuples} Given a commuting $d$ - tuple of bounded linear operators $\boldsymbol T= (T_1\ldots , T_d)$, our main goal is to find what natural conditions must be imposed on $\boldsymbol T$ to ensure that the commutators $[T_j^* , T_i]$, $1\leq i,j \leq d$,  are in trace class. Several conditions on $\boldsymbol T$ are known (see \cite{Atjoint, Douglasyan}) which achieve this goal. The theorem below provides an easy construction of such commuting $d$ - tuple of operators starting from a single operator $T$. First, we make a simple observation.

\begin{prop} \label{prop 2.1} 
Let $\boldsymbol{A}=(A_{1},\ldots, A_d)$ be a commuting projectively hyponormal (resp. hyponormal) tuple of operators on a Hilbert space $\mathcal{H}_{1}$ and $T$ be a hyponormal operator on (possibly) some other Hilbert space $\mathcal{H}_{2}$. Then $\boldsymbol{A}\otimes T=(A_{1}\otimes T ,\ldots, A_{d}\otimes T)$ is a commuting projectively  hyponormal (resp. hyponormal) tuple of operators on $\mathcal{H}_{1}\otimes \mathcal{H}_{2}$.
\end{prop} 
\begin{proof} For the first half of the proof, assume that $\boldsymbol A$ is projectively hyponormal.  
  For an arbitrary choice of $\alpha= (\alpha_1, \ldots, \alpha_d)\in \mathbb{C}^d$, let $\boldsymbol{A}_\alpha:= \sum_{i=1}^d \alpha_i A_i$ and $(\boldsymbol{A}\otimes T)_\alpha:= \sum_{i=1}^d \alpha_i 
  (A_i\otimes T) = \boldsymbol{A}_\alpha  \otimes T$. 
  Since $\boldsymbol{A}$ is commuting projectively hyponormal, we have 
  $[(\boldsymbol{A}_\alpha)^*, (\boldsymbol{A}_\alpha)]\geq 0$ for each $\alpha \in \mathbb{C}^d$. Now, 
  \begin{align*}
      \big [(\boldsymbol{A}_\alpha \otimes T)^*, (\boldsymbol{A}_\alpha\otimes T)\big ]& = (\boldsymbol{A}_\alpha)^*(\boldsymbol{A}_\alpha) \otimes T^*T-(\boldsymbol{A}_\alpha)(\boldsymbol{A}_\alpha)^* \otimes T T^*\\
    & \geq  (\boldsymbol{A}_\alpha)^*(\boldsymbol{A}_\alpha) \otimes T T^*-(\boldsymbol{A}_\alpha)(\boldsymbol{A}_\alpha)^* \otimes T T^* \\
   & = [(\boldsymbol{A}_\alpha )^*, (\boldsymbol{A}_\alpha)]\otimes T T^*\\
   & \geq 0.
  \end{align*}
  A similar proof: $\big[\!\!\big [(\boldsymbol A\otimes T)^*,\boldsymbol A\otimes T \big ]\!\!\big]\geq \big[\!\!\big[\boldsymbol A^*, \boldsymbol A\big ]\!\!\big]\otimes T T^*\geq 0$ completing the proof in the case $\boldsymbol A$ is hyponormal.
\end{proof}
What additional condition must be imposed on a hyponormal $d$ - tuple $\boldsymbol A\otimes T$ to conclude that 
$$\big[\!\!\big [(\boldsymbol A\otimes T)^*,\boldsymbol A\otimes T \big ]\!\!\big]: (\mathcal H_1\otimes \mathcal H_2)\otimes \mathbb C^d \to (\mathcal H_1\otimes \mathcal H_2)\otimes \mathbb C^d $$ is in trace class? The following theorem provides a set of very restrictive conditions ensuring this. 

\begin{thm}\label{Simpletensor}
Let $\boldsymbol A:=(A_{1},\ldots, A_{d})$ be a $d$ - tuple of commuting normal Hilbert-Schmidt operators on a Hilbert space $\mathcal{H}_{1}$ and $T$ be a hyponormal operator on a second (possibly different) Hilbert space $\mathcal{H}_{2}$. Assume that $T$ is $m$- polynomially cyclic.
Then $\boldsymbol A\otimes T:=(A_{1}\otimes T,\ldots , A_{d}\otimes T)$ is a commuting hyponormal $d$- tuple of operators and 
$$\big \|\big [A_{j}^{*}\otimes T^{*}, A_{i}\otimes T\big]\big \|_1 \leq \tfrac{m}{\pi} \nu \big(\sigma(T)\big)\|A_i\|_2 \|A_j\|_2, 1\leq i,j \leq d,$$
where $\nu$ is the normalized Lebesgue measure.
\end{thm}

\begin{proof}
 Since $T$ is a hyponormal and $m$- polynomially cyclic operator  on $\mathcal{H}_{2}$, by the Berger-Shaw theorem,  $\text{trace}\,([T^{*},T])\leq \tfrac{m}{\pi}\nu(\sigma(T))$. For $1\leq i,j \leq d,$ using Fuglede theorem for commuting normal operators, we have
	\begin{align} \label{Fuglede}
	\big [(A_{j}\otimes T)^*,(A_{i}\otimes T)\big ]=A_{j}^*A_{i}\otimes [T^*,T]=A_{i}A_{j}^*\otimes [T^*,T] 
	\end{align}
	Using Equation \eqref{Fuglede}, we obtain the equality
	\begin{align} \label{eqn:2.4}
\big[\!\!\big[ (\boldsymbol A\otimes T)^*, \boldsymbol A\otimes T \big]\!\!\big]
&=\begin{pmatrix}
  A_1 \\
  \vdots\\
  A_d
\end{pmatrix} \begin{pmatrix}
  A_1^*&\ldots & A_d^*
\end{pmatrix} \otimes[T^*,T].
\end{align}
Since  $T$ is hyponormal by hypothesis, it follows that the commuting tuple $\boldsymbol A \otimes T$ is hyponormal.

Also, since $A_i$, $1\leq i \leq d$, is a Hilbert-Schmidt operator, it follows that $A_i A_j^*$, $1\leq i, j \leq d$, is in trace-class. The operator $[T^*, T]$ is in trace-class,  consequently (using Equation \eqref{Fuglede}), we conclude that 
$\big [A_{j}^{*}\otimes T^{*}, A_{i}\otimes T\big ]$, $1 \leq i,j \leq d,$ is in trace-class. Moreover, we have 
\begin{align*}
    \big \|\big [(A_{j}\otimes T)^*,(A_{i}\otimes T)\big ]\big \|_1&=\big \|A_{j}^*A_{i}\otimes \big [T^*,T\big ]\big \|_1\\& = \text{\rm trace} \big ( | A_j^* A_i | \big ) \text{\rm trace} \big [T^*,T\big ] \\
    & \leq \tfrac{m}{\pi} \nu \big(\sigma(T)\big)\big \|A_i\big \|_2 \big \|A_j\big \|_2,
\end{align*}
where the last inequality follows since $\|AB\|_1 \leq \|A\|_2 \|B\|_2$ for any pair of Hilbert - Schmidt operators $A$ and $B$. 
\end{proof}
The following lemma is certainly well known, however for the shake of completeness, we provide a proof.
\begin{lem}\label{lem:1.3}
Let $X,Z$ be positive operators on some Hilbert space $\mathcal H_1$ and $Y$ be a nonzero operator on (possibly) some other Hilbert space $\mathcal{H}_2$. Assume that either 
$$X\otimes Y^* Y \geq Z \otimes Y Y^*, \mbox{\rm ~or ~} Y^* Y \otimes X  \geq  Y Y^*\otimes Z.$$
Then $X \geq Z$.
\end{lem}
\begin{proof}
    Suppose $X\otimes Y^* Y \geq Z \otimes Y Y^*$. For any pair  of unit vectors $u\in \mathcal{H}_1$ and $v\in \mathcal{H}_2$, we have 
    \begin{align}
       \langle X u , u \rangle \|Y v\|^2 &= \langle X\otimes Y^*Y u\otimes v , u\otimes v\rangle \nonumber \\ &\geq \langle Z \otimes Y Y^* u\otimes v , u\otimes v\rangle =  \langle Z u , u \rangle \|Y^* v\|^2.
        \end{align}
        Taking supremum over all unit vectors $v\in \mathcal H_2$, we see that $X \geq Z$. The proof is similar if we start with the inequality $Y^* Y \otimes X  \geq  Y Y^*\otimes Z.$
        \end{proof}
In passing, we note that there is a converse to Proposition \ref{prop 2.1}, which we give below. 
\begin{prop}\label{prop:1.4}
Let $\boldsymbol A= (A_1, \ldots, A_d)$ be $d$- tuple of operators on a Hilbert space $\mathcal{H}_1$ and $T$ be a nonzero operator on (possibly) some other Hilbert space $\mathcal{H}_2$. Assume that the tuple of operators $\boldsymbol A\otimes T:= ( A_1 \otimes T, \ldots , A_d \otimes T)$ on $\mathcal{H}_1\otimes \mathcal{H}_2$ is  commuting and projectively  hyponormal (resp. hyponormal).  Then we have the following.
\begin{enumerate}
\item[(i)] Each of the operators $A_1, \ldots, A_d$ is hyponormal on $\mathcal{H}_1$ and the operator $T$ is hyponormal on  $\mathcal{H}_2$.
\item[(ii)]  The tuple of operators $\boldsymbol{A}$ is  commuting and projectively hyponormal (resp. hyponormal) on $\mathcal H_1$.
\end{enumerate}
\end{prop}
\begin{proof} Since $\boldsymbol A\otimes T$ is projectively hyponormal, the commutators 
\begin{align*}
[A_i^*\otimes T^*, A_i\otimes T]  &= (A_i^{*}A_i\otimes T^{*}T - A_i A_i^{*}\otimes T T^{*})\geq 0.  
\end{align*}
First, applying Lemma \ref{lem:1.3} with $X=A_i^*A_i$, $Z=A_i A_i^*$ and $Y=T$, and then with $X=T^*T, Z=T T^*$ and $Y=A_i$, we conclude that  $A_i$, $1\leq i \leq d$, and $T$ are hyponormal. By hypothesis, 
\begin{align*}
0 &= [(A_i\otimes T),(A_j\otimes T)]\\
&= [A_i, A_j] \otimes T^2
\end{align*}
whenever $i\neq j$. This implies: either $[A_i,A_j]$=0 or $T^{2}=0$. If $T^{2}=0$, then the spectral radius of $T$, $r(T)=0$. Since $T$ is a hyponormal operator we have $r(T) = \|T\|$ (see \cite[Theorem 1]{stampfli}). So $\|T\|=0$ which is a contradiction. Hence we have $[A_i,A_j]$=0.

For an arbitrary choice of complex numbers $\alpha_1, \ldots, \alpha_d$, let $\boldsymbol{A}_\alpha:= \sum_{i=1}^d \alpha_i A_i$. If $\boldsymbol A\otimes T$ is assumed to be projectively hyponormal, then $[(\boldsymbol A_{\alpha}\otimes T)^*, \boldsymbol A_{\alpha}\otimes T]\geq 0.$ Now, setting  
$$X := \boldsymbol{A}_\alpha^* \boldsymbol{A}_\alpha,\,\, Z:= \boldsymbol{A}_\alpha \boldsymbol{A}_\alpha^*\,\, \text{and}\,\, Y:= T,$$
and applying Lemma \ref{lem:1.3}, it follows that $\boldsymbol{A}$ is projectively hyponormal.

A similar proof: Set $$X:= \begin{pmatrix}
  A_1^* A_1 & \ldots & A_d^* A_1\\
   \vdots    & \ddots   &  \vdots \\
    A_1^* A_d & \ldots & A_d^* A_d
\end{pmatrix}, \,\,\,\, Z:= \begin{pmatrix}
   A_1 A_1^* & \ldots & A_1 A_d^*\\
   \vdots    & \ddots   &  \vdots \\
    A_d A_1^* & \ldots & A_d A_d^*
\end{pmatrix},$$ both acting on the Hilbert space $\mathcal H_1\otimes \mathbb C^d$, and $Y:=T$. The hypotheses of Lemma \ref{lem:1.3} are met when $\boldsymbol A\otimes T$ is assumed to be hyponormal (see Equation \eqref{eqn:2.4}). Thus, in this case,  we have $X\geq Z$ 
completing the proof of the hyponormality of the commuting $d$ - tuple $\boldsymbol A$. 
\end{proof}

\subsection{Hyponormality  and trace class commutators}
For many of the $d$- tuple of operators discussed in the previous section, we show that the conclusion of the Douglas-Yan theorem, which is recalled below, remains valid even though, they don't necessarily satisfy the hypotheses of the theorem. One of these hypotheses is that of the Krull dimension of the quotient $\mathbb C[z_1, \ldots , z_d]/\mathcal I_{\boldsymbol T}$, where $$\mathcal I_{\boldsymbol T}= \big \{p \in \mathbb C[z_1, \ldots , z_d]:p(\boldsymbol T)=0 \big \}$$ is the \emph{ vanishing ideal} of the $d$ - tuple $\boldsymbol T$. 
\begin{defn}[Krull dimension] 
The Krull dimension of a ring $\mathcal R$ is defined to be the maximum of those positive integers $n$ for which there is an ascending chain of prime ideals of the form  
$$\{0\}= P_{0}\subsetneq P_1 \subsetneq \ldots \subsetneq P_n \subsetneq \mathcal R.$$
\end{defn}
We also give some natural examples of commuting $d$- tuples of operators satisfying the hypothesis of the Douglas-Yan theorem. Finally, we indicate a proof of this theorem with the weaker hypothesis of projective hyponormality.  

\begin{thm}
Let $\boldsymbol T=(T_1, \ldots, T_d)$ be a hyponormal $d$- tuple of operators on the Hilbert space $\mathcal{H}$ such that $\boldsymbol T$ is finitely polynomially cyclic. 
If the Krull dimension of $\mathbb C[z_1, \dots ,z_d]/\mathcal I_{\boldsymbol T}$ is $1$, then $[T_j^*, T_i]$ is in trace-class for all  $i,j$$(1\leq i,j \leq d)$.  
\end{thm}

An algebraic closed set $V\subseteq \mathbb C^d$ is
the  set of common zeros of a set of polynomials in $\mathbb C[z_1, \dots ,z_d]$. Furthermore,  if $\mathcal I_V \subseteq \mathbb C[z_1, \dots ,z_d]$ is the ideal of all polynomials vanishing on the set $V$, then the set $V$ is called an algebraic curve if $\dim \big ( \mathbb C[z_1, \dots ,z_d] / \mathcal I_V\big)=1$. A closed algebraic set $V$ is said to be thin if it is a subset of an algebraic curve. 
In the special case, when the $d$- tuple $\boldsymbol T$ is subnormal, $\mbox{dim}\big ( \mathbb{C}[z_1,\ldots, z_d]/\mathcal I_{\boldsymbol T}\big )=1$ is equivalent to the fact that the Taylor joint spectrum of the $d$- tuple is thin, see \cite{Douglasyan}.

It is proved in Theorem \ref{Simpletensor} that the commutators $[( A_j \otimes T)^*, A_i \otimes T]$, $1\leq i,j \leq d$, 
belong to the trace class after having imposed suitable hypothesis on the commuting tuple $\boldsymbol A$ and the operator $T$.  This set of hypotheses does not include polynomial cyclicity of the commuting $d$- tuple $\boldsymbol A \otimes T$. Polynomial cyclicity is one of the two main assumptions in the Douglas-Yan theorem. However, we show that the Taylor joint spectrum of the commuting $d$- tuple $\boldsymbol A \otimes T$ appearing in Theorem \ref{Simpletensor} is a union of thin sets. This was the other hypothesis in the Douglas-Yan Theorem. The Theorem stated below is slightly more general than \cite[Proposition 3.16]{BGMS}, however, the proof is similar and therefore omitted. 
 \begin{prop}\label{dspectrum}
Let $\big \{\boldsymbol S_i:=\big (S_{i1}, \ldots , S_{id}\big )\big \}_{i=1}^n$, where $n\in \mathbb N\cup \{\infty\}$, be a set of commuting $d$- tuple of bounded linear operators on the Hilbert spaces $\mathcal{H}_{i}$. Then the Taylor joint spectrum  $\sigma(\oplus_{i=1}^n \boldsymbol S_i)$ of the direct sum $\oplus_{i=1}^n \boldsymbol S_i$ equals $\overline{\cup_{i=1}^n\sigma(\boldsymbol S_i)}$. 
\end{prop}
The Theorem below follows directly from Proposition \ref{dspectrum}. 
\begin{thm} \label{DYexp}
Let $\boldsymbol A:=(A_{1},\ldots, A_{d})$ be a $d$- tuple of commuting normal Hilbert-Schmidt operators on a Hilbert space $\mathcal{H}_{1}$ and $T$ be an operator on a second (possibly different) Hilbert space $\mathcal{H}_{2}$.
Then the joint spectrum of the commuting $d$- tuple of operators $\boldsymbol A\otimes T:=(A_{1}\otimes T,\ldots , A_{d}\otimes T)$ is a countable union of thin sets. 
\end{thm}
\begin{proof}
Let $\boldsymbol \lambda^{(k)} = (\lambda^{(k)}_i)_{i\in I}$, $1\leq k \leq d$, $I$ countable, be the set of eigenvalues of $A_k$.
  We can assume, without loss of generality, that $A_k= \sum \lambda^{(k)}_i P^{(k)}_i$, where $P^{(k)}_i$ is the orthogonal projection to the eigenspace corresponding to $\lambda^{(k)}_i$. Consequently, the tuple $\boldsymbol A\otimes T$ is the direct sum 
$$ \bigoplus_{i\in I} \big ( \lambda^{(1)}_i  T\,, \ldots ,\, \lambda^{(d)}_i  T \big ). $$
For a fixed $i\in I$, the spectrum of the operator $\big ( \lambda^{(1)}_i  T\,, \ldots ,\, \lambda^{(d)}_i  T \big )$ is the set $\{ (\lambda^{(1)}_i z , \ldots , \lambda^{(d)}_i  z): z\in  \sigma(T) \} \subseteq \mathbb C^d$.  Now, applying  Proposition \ref{dspectrum}, we conclude that the spectrum of the $d$- tuple of operators $\boldsymbol A\otimes T$ is the set 
$$\sigma(\boldsymbol A\otimes T)= \bigcup_{i\in I} \big \{ (\lambda^{(1)}_i z , \ldots , \lambda^{(d)}_i  z): z\in  \sigma(T) \big \}.$$
Thus $\sigma(\boldsymbol A\otimes T)$ is a union of at most countably many thin sets.
\end{proof}
\begin{rem} The simple examples of Theorem \ref{Simpletensor} show that the conclusion of membership in the trace class can be achieved without imposing either of the hypotheses  mandated in the Douglas-Yan theorem. 

 As shown in Theorem \ref{DYexp}, the spectrum of the $d$ - tuple of commuting operators $\boldsymbol A\otimes T$, in general, is only a countable union of thin sets. Although, if $\boldsymbol A$ is a commuting $d$- tuple with the property:
 $$\mbox{\rm card} \big \{ \boldsymbol \lambda := \big ( \lambda^{(1)}  , \ldots ,\, \lambda^{(d)} \big ) \mid \boldsymbol \lambda \in \sigma(\boldsymbol A) \big \} < \infty, $$
 then the $\sigma(\boldsymbol A\otimes T)$ is a union of finitely many thin sets and therefore it is thin. 
 \end{rem}
 The following proposition gives a different set of examples satisfying the hypotheses of the Douglas-Yan theorem.  
 
 Recall that an operator $T$ (resp. $d$ - tuple) of hyponormal operator is said to be \textit{pure} if it has no reducing subspace $\mathcal{H}_0$ such that $T$ restricted to $\mathcal{H}_0$ is normal.
 
\begin{prop}
If $T$ is a pure hyponormal operator, and $\mathcal I_{\boldsymbol T} \subseteq \mathbb C[\boldsymbol z]$ is the vanishing ideal of $\boldsymbol T:=(T,\ldots, T^d)$, then $\dim\big ( \mathbb C[\boldsymbol z]/{\mathcal I_{\boldsymbol T}} \big ) = 1$. Moreover, if $T$ is $m$ - cyclic, then  $(T,\ldots, T^d)$ is $n$ - cyclic for some $n$, $n \leq m$. 
\end{prop}
\begin{proof} Clearly, if $T$ is $m$- polynomially cyclic, then $(T,\ldots ,T^d)$ is $n$- polynomially cyclic for some $n$ with $n \leq m$.
By definition the ideal $\mathcal I$ generated by the polynomials $p^{(1)}_T,\ldots, p^{(d-1)}_T: p^{(1)}_T(\boldsymbol z)= z_1^2 - z_2,\ldots, p^{(d-1)}_T(\boldsymbol z)= z_1^d - z_d$ is included in $\mathcal I_{\boldsymbol T}$. 

The kernel of the map $\Phi:\mathbb{C}[z_1,\ldots, z_d]\to \mathbb{C}[z]$ given by $(\Phi\, f)(z)= f(z, z^2,\ldots, z^d)$ is $\mathcal I$ and therefore, $\mathbb{C}[z_1,\ldots, z_d]/\mathcal I \cong \mathbb{C}[z]$. It follows that the Krull dimension of $\mathbb{C}[z_1,\ldots, z_d]\big / \mathcal I$ equals that of $\mathbb{C}[z]$, which is $1$. Therefore, to complete the proof, all we have to do is show that  $\mathcal I_{\boldsymbol T}$ is included in $\mathcal I$. 

If there is a polynomial $p$ not in $\mathcal I$ such that $p(T, T^2,\ldots, T^d)= 0$, then the polynomial $\phi(p)(z) := p(z,z^2, \ldots, z^d)$ is not the zero polynomial and $\phi(p)(T)=0$. But $\phi(p)(T)=0$ if and only if the spectrum of $T$ is finite. Now, if  $T$ is a pure hyponormal operator, then the spectrum of $T$ can not be discrete hence cannot be finite (see \cite[Cor. 2]{stampfli}). Therefore, the vanishing ideal of the $d$- tuple $(T,\ldots, T^d)$ is the ideal generated by the polynomials ${z_1}^2-z_2, \ldots, {z_1}^d-z_d$.  This completes the proof. 
\end{proof}

Recall that Athavale, in the paper \cite{Atjoint}, had introduced the notion of the  projectively hyponormal operators, although, he called them weakly hyponormal. Among other things, the proof of the inequality \eqref{eqn:2.6} appears in the proof of Proposition 4 in his paper. However, the estimate for $|\text{Trace}[T_2^*,T_1]|$ that follows is valid if $[T_2^*,T_1]$ is trace class. This was claimed in the proposition but not verified. We provide this proof here. For our proof, it is necessary to adopt the following slightly different definition of the trace norm $\|T\|_1$, see \cite[Proposition 3.6.5]{Simon}, which however is equivalent to the expression in Definition \ref{defn1.1}:
\begin{equation}\label{tracenorm} 
\|T\|_1 = \sup_{\{f_n\},\{g_n\}} \sum_{n=1}^\infty |\langle T f_n , g_n \rangle|, 
\end{equation}
where the supremum is taken over all pairs $\{f_n\}$ and $\{g_n\}$ of orthonormal sets in the Hilbert space $\mathcal H$.

\begin{lem} \label{Estimate for A} Let $\boldsymbol T=(T_1,\ldots,T_d)$ be a commuting projectively hyponormal $d$- tuple of operators on the Hilbert space $\mathcal{H}$. Furthermore, assume that for $1\leq k \leq d$, each $[T_k^*,T_k]$ is in trace-class. Then $[T_j^*,T_k]$ is also in trace-class for all $1\leq k,j\leq d$. 
\end{lem} 
\begin{proof}  From the projective hyponormality of $\boldsymbol T$, for any fixed pair of indices $j,k$, it follows that 
\begin{align}\label{15}|\alpha_j|^2\langle f,[T_j^*,T_j]f\rangle+2\text{Re}\,\bar{\alpha_j}\alpha_k\langle f,[T_j^*,T_k]f\rangle+ |\alpha_k|^2 \langle f,[T_k^*,T_k]f\rangle \geq 0, \alpha_j, \alpha_k \in \mathbb{C},  f\in \mathcal{H}, \end{align}
where $\alpha_j, \alpha_k$ are an arbitrary pair of complex numbers. 
By a reasoning identical to the derivation of the Cauchy-Schwarz inequality, we get 
	\begin{align}\label{eqn:2.6}
	|\langle f,[T_j^*,T_k]f\rangle|^2\leq \langle f,[T_j^*,T_j]f\rangle \langle f,[T_k^*,T_k]f\rangle, \,f\in \mathcal{H}. 
	\end{align}
For $f,g\in \mathcal H$, using the polarization identity and Equation \eqref{eqn:2.6} together, we have 
	\begin{align*}
	4|\langle f,[T_j^*, T_k]g\rangle|&\leq |\langle f+g, [T_j^*, T_k](f+g)\rangle|+|\langle f-g,[T_j^*, T_k](f-g)\rangle| \nonumber \\&\phantom{Param}+|\langle f+ig, [T_j^*, T_k](f+ig)\rangle|+|\langle f-ig, [T_j^*, T_k](f-ig)\rangle|\\
	& \leq (\langle f+g,[T_k^*,T_k](f+g)\rangle \langle f+g,[T_j^*,T_j](f+g)\rangle) ^\frac{1}{2} \\
	& \phantom{Param} + (\langle f-g,[T_k^*,T_k](f-g)\rangle \langle f-g,[T_j^*,T_j](f-g)\rangle) ^\frac{1}{2} \\
& \phantom{Paramita} + (\langle f+ig,[T_k^*,T_k](f+ig)\rangle \langle f+ig,[T_j^*,T_j](f+ig)\rangle) ^\frac{1}{2}\\
	& \phantom{ParamitaM}+ (\langle f-ig,[T_k^*,T_k](f-ig)\rangle \langle f-ig,[T_j^*,T_j](f-ig)\rangle) ^\frac{1}{2}.	
	\end{align*}
	Next using the inequality $xy  \leq \tfrac{x^2 + y^2}{2}$, $x,y\in \mathbb R_{+}$ and collecting the appropriate terms back for another application of polarization identity we get that
	\begin{align*}
	4|\langle f,[T_j^*, &T_k]g \rangle| \leq \frac{1}{2} \big(\langle f+g,[T_k^*,T_k](f+g)\rangle +\langle f+g,[T_j^*,T_j](f+g)\rangle + \langle f-g,[T_k^*,T_k](f-g)\rangle 
	\\& \phantom{Par}+\langle f-g,[T_j^*,T_j](f-g)\rangle + \langle f+ig,[T_k^*,T_k](f+ig)\rangle +\langle f+ig,[T_j^*,T_j](f+ig)\rangle 
	\\&\phantom{Param}	+\langle f-ig,[T_k^*,T_k](f-ig)\rangle +\langle f-ig,[T_j^*,T_j](f-ig)\rangle\big) 
	\\& = {2} \big (\langle f,[T_k^*,T_k]f\rangle+\langle g,[T_k^*,T_k]g\rangle+\langle f,[T_j^*,T_j]f\rangle+\langle g,[T_j^*,T_j]g\rangle \big ),	
	\end{align*} where we have also used the fact that $[T_j^*, T_j]\geq 0$ for all $j$.
	 For any pair of orthonormal sets $\{f_n\}$ and $\{g_n\}$ ,
    \begin{align*}
	\sum_n \big |\langle f_n, [T_j^*, T_k]g_n\rangle \big |
	 \leq & \frac{1}{2} \big (\sum_{n}\langle f_n, [T_k^*, T_k]f_n\rangle+\sum_{n}\langle g_n, [T_k^*,T_k]g_n\rangle
	\\  &\phantom{Param} +\sum_{n} \langle f_n, [T_j^*,T_j]f_n\rangle+\sum_{n} \langle g_n, [T_j^*, T_j]g_n\rangle \big)\\
	\leq & \frac{1}{2}\big (2\big \|[T_k^*, T_k]\big \|_1+2\big \|[T_j^*, T_j]\big \|_1\big )\\
	= & \big \|[T_k^*, T_k]\big \|_1+\big \|[T_j^*, T_j]\big \|_1.\end{align*}
	Taking supremum over all possible orthonormal sets and using \eqref{tracenorm} we get that
	\begin{align}\label{eqn:2.7}
	\big \|[T_j^*, T_k]\big \|_1=\big \|[T_k^*, T_j]\big \|_1 &=\sup_{\{f_n\}\{g_n\}}\sum_{n} |\langle [T_k^*, T_j] f_n, g_n\rangle|\nonumber\\&=
	\sup_{\{f_n\}\{g_n\}}\sum_{n} |\langle f_n, [T_j^*, T_k]g_n\rangle| \leq \big \|[T_k^*, T_k]\big \|_1+\big \|[T_j^*, T_j]\big \|_1 < \infty.\end{align}
	Thus $ [T_j^*, T_k]$	is in trace class.		
\end{proof}

For $i\not = j$, the estimate of the trace norm $[T_j^*,T_i]$, $1\leq i,j \leq d$, obtained in Theorem \ref{Asecond} is less than or equal to half of the estimate of the same norm obtained in Theorem \ref{Afirst}, which is only natural because of the stronger hypotheses in Theorem \ref{Asecond}.


\begin{thm} \label{Afirst}
    Let $\boldsymbol T=(T_1, \ldots , T_d)$ be a projectively hyponormal $d$- tuple of  operators. Suppose that each $T_i$, $1\leq i \leq d$, is $m_i$- polynomially cyclic. Then the operators $[T_j^*,T_i]$ are in trace class. Moreover, we have 
    $$\big \|[T_j^*,T_i]\big \|_1 \leq  \begin{cases} \tfrac{m_i}{\pi}\nu(\sigma(T_i)) &  \mbox{\rm ~if~} i = j \\ \tfrac{m_i}{\pi} \nu(\sigma(T_i)) +  \tfrac{m_j}{\pi} \nu(\sigma(T_j))& \mbox{\rm ~if~} i\not = j, \end{cases} $$
    where $\nu(\sigma(T_i))$  is the Lebesgue measure of the spectrum of $T_i$. 
\end{thm}
The verification of the estimate for the trace norm of $[T_j^*,T_i]$ in Theorem \ref{Afirst} assuming that  $(T_1, \ldots T_d)$ is projectively hyponormal follows directly from Equation \eqref{eqn:2.7} and the Berger-Shaw theorem. 

\begin{thm}\label{Asecond}
    Let $\boldsymbol T=(T_1, \ldots , T_d)$ be a hyponormal $d$- tuple of  operators. Suppose that each $T_i$, $1\leq i \leq d$, is $m_i$- polynomially cyclic. Then the operators $[T_j^*,T_i]$ are in trace class. Moreover, we have 
    $$\big\|[T_j^*,T_i]\big \|_1 \leq  \begin{cases} \tfrac{m_i}{\pi}\nu(\sigma(T_i)) &  \mbox{\rm ~if~} i = j \\ \big (\tfrac{m_i}{\pi} \nu(\sigma(T_i))\big )^{\frac{1}{2}}   \big ( \tfrac{m_j}{\pi} \nu(\sigma(T_j))\big )^{\frac{1}{2}}& \mbox{\rm ~if~} i\not = j, \end{cases} $$
    where $\nu \big(\sigma(T)\big)$  is the Lebesgue measure of the spectrum of $T$.  
\end{thm}
\begin{proof}
    Hyponormality of the $d$- tuple $\boldsymbol T$ gives, for any pair $x, y \in \mathcal{H}$,
    $$\big |\langle [T_j^*, T_i] x, y\rangle \big |\leq {\langle [T_i^*, T_i] x, x\rangle}^{\frac{1}{2}} {\langle [T_j^*, T_j] y, y\rangle}^{\frac{1}{2}}.$$
    For any pair of orthonormal sets $\{f_n\}$ and $\{g_n\}$ ,
    \begin{align*}
	\sum_n \big |\langle [T_j^*, T_i]f_n, g_n\rangle \big |&\leq \sum_n{\langle [T_i^*, T_i] f_n, f_n\rangle}^{\frac{1}{2}} {\langle [T_j^*, T_j] g_n, g_n\rangle}^{\frac{1}{2}}\\&\leq \big(\sum_n\langle [T_i^*, T_i] f_n, f_n\rangle \big)^{\frac{1}{2}} \big(\sum_n\langle [T_j^*, T_j] g_n, g_n\rangle \big)^{\frac{1}{2}}\\& \leq \big\|[T_i^*, T_i]\big\|_1^{\frac{1}{2}}\big\|[T_j^*, T_j]\big\|_1^{\frac{1}{2}}.
    \end{align*}
    Now the conclusion follows from the definition of the trace norm given in \eqref{tracenorm} and the Berger-Shaw theorem. 
\end{proof}

We point out that Lemma \ref{Estimate for A} appears in \cite{Douglasyan} with the stronger assumption of hyponormality of the $d$ - tuple $\boldsymbol T$. Apart from Lemma \ref{Estimate for A}, we need a second preparatory lemma to prove the Douglas-Yan theorem with the weaker hypothesis of projective hyponormality. 
\begin{lem} \label{Estimate for B} Let $\boldsymbol T=(T_1,\ldots, T_d)$ be a projectively hyponormal commuting  $d$- tuple of operators on a Hilbert space $\mathcal{H}$ and  $a_{i j}$, $1\leq i \leq k$, $1\leq j \leq d$;  be a set of $k d$ complex scalars. Suppose $g_i(\boldsymbol T)=\sum_{j=1}^{d}a_{i j}T_j,\, i=1,\ldots,k$. Then $\big (g_1(\boldsymbol T),\ldots, g_k(\boldsymbol T)\big )$ is also projectively hyponormal on $\mathcal{H}$.	
\end{lem}
\begin{proof} $\boldsymbol T$ is projectively hyponormal implies the operator   $\alpha_1T_1+\alpha_2T_2+\ldots+\alpha_d T_d$ is hyponormal for all $\alpha
	_i\in\mathbb{C}$. Since  $\sum_{i=1}^{k}\lambda_i g_i(\boldsymbol T)=\sum_{j=1}^{d}(\sum_{i}^k\lambda_i a_{i j}) T_j,$ therefore $\big (g_1(\boldsymbol T),\ldots,g_k(\boldsymbol T)\big )$ is also projectively hyponormal.
\end{proof}

We state the following theorem without proof, which is a restatement of \cite[ Theorem 2]{Douglasyan} with the weaker hypothesis of projective hyponormality instead of hyponormality. 
The proof is along the same lines as the original proof 
combined with Lemma \ref{Estimate for A} and  Lemma \ref{Estimate for B} and therefore 
omitted. 
\begin{thm}\label{12345}
	Let $\boldsymbol T=(T_1,\ldots,T_d)$ be a  projectively hyponormal commuting $d$- tuple of operators on a Hilbert space $\mathcal{H}$ such that $\boldsymbol T$ is $m$-  cyclic. Assume that $\mathbb{C}[z_1,\ldots,z_d]/\mathcal I_{\boldsymbol T}$ has Krull dimension 1, where $\mathcal I_{\boldsymbol T}$ is the vanishing ideal of $\boldsymbol T$, then $[T_j^*,T_i]$ is in trace class for all $1\leq i,j\leq d$.
\end{thm}

However note that finite cyclicity of $\boldsymbol T$ has no implication for the finite cyclicity of the individual operators $T_i$, $1\leq i \leq d$. Consequently, unlike in Theorem \ref{Afirst}, there is no estimate of the trace norm of $[T^*_j,T_i]$ in Theorem \ref{12345}.

In the proof of the following theorem, we consider those polynomials $p$ in $\mathbb C[z_1, \ldots , z_d]$ with the property that none of the partial derivatives $\partial_k p$, $1\leq k \leq d$, is the zero polynomial. We call such a polynomial \textit{full}. 

Suppose $(T_1,T_2)$ is any pair of commuting operators and $p(T_1,T_2) = 0$. If $p$ is not full, then  $p$ is of the form: (a) $p(z_1,z_2) = \sum_{k=1}^m a_k z_1^k$, or (b) 
$p(z_1,z_2) = \sum_{k=1}^n b_k z_2^k$. In either case, the spectrum of $T_1$ or $T_2$ is finite. Now, if  $T_1$ and $T_2$ are pure
hyponormal operators, then the spectrum of neither of these can be discrete hence cannot be finite (see \cite[Cor. 2]{stampfli}). This contradiction shows that for any commuting pair $T_1$, $T_2$ of pure hyponormal operators, if $p$ is not a full polynomial, then $p(T_1,T_2)$ cannot be zero. 
\begin{thm} \label{strengthening}
Let $\boldsymbol T=(T_1, T_2)$ be a (pure) projectively hyponormal pair of commuting operators on the Hilbert space $\mathcal{H}$ such that $\boldsymbol T$ is $m$- polynomially cyclic. Furthermore, assume that there exists a polynomial $p \in \mathbb{C}[z_1, z_2]$ such that $p(T_1, T_2)=0,$ then $[T_j^*,T_i]$ is in trace class for all $1\leq i,j\leq 2.$
\end{thm}
\begin{proof}The discussion preceding the theorem shows that if $p(T_1, T_2) = 0$
for a polynomial $p$, then it must be full. 
Let degree of $z_2$ in $p$ be $k$.
For any polynomial $q\in \mathbb{C}[z_1, z_2]$ using the Division algorithm in $(\mathbb C [z_1])[z_2]$  we get 
$$q(z_1, z_2)=q_1(z_1, z_2)p(z_1, z_2)+r(z_1, z_2)$$
where degree of $z_2$ in $r$ is less than $k$. 
Since $p(T_1, T_2)=0$ it follows that $q(T_1, T_2)= r(T_1, T_2).$ We conclude that $T_1$ is $m k$- polynomially cyclic.
Now projective hyponomality of the pair $\boldsymbol T$ implies  that $T_1$ is hyponormal. Hence, by the Berger-Shaw theorem, $[T_1^*, T_1]$ is in trace-class. Similarly, since $T_2$ is also polynomially cyclic, one can prove $[T_2^*, T_2]$ is in trace-class. Finally, from Lemma \ref{Estimate for A} we conclude that for $1\leq i,j \leq 2$, $[T_j^*, T_i]$ is in trace-class.
\end{proof}
It is possible to construct a large class of commuting pairs of operators, using Theorem \ref{12345}, where the Douglas-Yan theorem applies with the slightly weaker hypothesis of projective hyponormality.

\begin{cor}
	Let $\boldsymbol T=(T_1, T_2)$ be a  projectively hyponormal commuting pair of operators on a Hilbert space $\mathcal{H}$. Assume that $p(T_1,T_2) = 0$ for some $p\in \mathbb C[z_1,z_2]$. Then $\dim \big ( \mathbb C[z_1,z_2]/\mathcal I_{\boldsymbol T} \big ) = 1$, where $\mathcal I_{\boldsymbol T}$ is the vanishing ideal of $\boldsymbol T$.  
\end{cor}
\begin{proof}
  The only possibilities for the Krull dimension of $\mathbb C[z_1,z_2]/\mathcal I_{\boldsymbol T}$ are $1$ or $0$. But if the Krull dimension is $0$, then $\mathcal I_{\boldsymbol T}$ must be a maximal ideal. Since $T_i,\, i=1,2,$ is a pure hyponormal operator, it follows that no maximal ideal can be the vanishing ideal of $(T_1, T_2)$. Hence the $\text{dim}\,\mathbb C[z_1,z_2] /\mathcal I_{\boldsymbol T}=1.$
\end{proof}

\section{Determinant and the generalized commutator}\label{weighted}
Let $\mathbb B_d:=\{\boldsymbol z=(z_1, \ldots ,z_d)\in \mathbb C^d: |z_1|^2 +\cdots +|z_d|^2 < 1\}$ be the Euclidean ball. The set of monomials 
$$\big \{\tfrac{(\lambda)_{|\boldsymbol \alpha|}}{\boldsymbol \alpha!} z_1^{\alpha_1}\cdots  z_d^{\alpha_d}: \boldsymbol \alpha = (\alpha_1, \ldots ,\alpha_d) \in \mathbb N_0^d\big \},$$
is a linearly independent set of vectors in $\mathbb C[\boldsymbol z],$ where for $\lambda >0$, $(\lambda)_n := \lambda (\lambda+1) \cdots (\lambda +n-1)$ is the Pochhammer symbol.  The series 
$K^{(\lambda)}(\boldsymbol z,\boldsymbol w):= \sum_{\boldsymbol \alpha}(\lambda)_{|\boldsymbol \alpha|} \boldsymbol z^{\boldsymbol \alpha} \bar{\boldsymbol w}^{\boldsymbol \alpha}$
converges absolutely, uniformly on compact subsets of $\mathbb B_d \times \mathbb B_d$.
Furthermore, $K^{(\lambda)} (\boldsymbol z,\boldsymbol w) = (1-\langle \boldsymbol z, \boldsymbol w\rangle)^{-\lambda}$ defines  positive definite kernel on $\mathbb B_d$ which is holomorphic in the first and anti-holomorphic in the second variable. These are the weighted Bergman kernels of the Euclidean ball $\mathbb B_d$. They determine a Hilbert space of holomorphic functions $\mathcal H^{(\lambda)}(\mathbb B_d),$ $\lambda >0$, where $K^{(\lambda)}$ serves as a \emph{reproducing kernel}, that is,
$$ \big \langle f, K^{(\lambda)}(\cdot, \boldsymbol w) \big \rangle = f(\boldsymbol w),\,f\in \mathcal H^{(\lambda)}(\mathbb B_d),\,\boldsymbol w\in \mathbb B_d.$$ 
Let $M_i: \mathcal H^{(\lambda)}(\mathbb B_d) \to \mathcal H^{(\lambda)}(\mathbb B_d)$,  be the operator of multiplication by the coordinate function $z_i$, $1\leq i \leq d$.
Let $\boldsymbol M^{(\lambda)}$ denote the commuting $d$ - tuple $(M_1, \ldots , M_d)$.
It is not hard to verify that the commutators $[M_j^*,M_i]$, $1\leq i \leq d$, on the weighted Bergman spaces $\mathcal H^{(\lambda)}(\mathbb B_d)$, $\lambda \geq 1$, are in Schatten $p$ - class if and only if $p > d$. A related question involving compresion of these coordinate multiplications to the orthocomplement of an invariant subspace of the form $[\mathcal I]\subseteq \mathcal H^{(\lambda)}(\mathbb B_d)$ remains unanswered. (Here $[\mathcal I]$ denotes the closure in $\mathcal H^{(\lambda)}(\mathbb B_d)$ of the ideal $\mathcal I$.) Indeed, the Arveson-Douglas conjecture, see \cite[Conjecture 3]{GY}, addressing such questions is an active area of research. 
\begin{conj}[Arveson-Douglas] 
Assume $\mathcal I$ is a homogeneous ideal of the polynomial ring $\mathbb C[z_1, \ldots, z_d]$ and $[\mathcal I]$ is the closure of $\mathcal I$ in $\mathcal H^{(\lambda)}(\mathbb B_d)$. Then for all $r > \dim_{\mathbb C} \mathcal Z(\mathcal I)$, where $\mathcal Z(\mathcal I)$ is the common zero set of $\mathcal I$, the quotient module $[\mathcal I]^\perp$ is r-essentially normal.
\end{conj}
The $d$ - tuple $\boldsymbol M^{(\lambda)}$ is known to be hyponormal if and only if $\lambda \geq d$, see \cite[pp. 605]{CY}. It is cyclic with cyclic vector $1$. Therefore, the obvious generalization of the Berger-Shaw theorem fails for the $d$ - tuple $\boldsymbol M^{(\lambda)}$.
 So, we look for a different set of criteria for a $d$ - tuple $\boldsymbol T$  of operators and a suitable function 
 $$\mathfrak D:\mathcal B\big (\mathcal H\otimes \mathbb C^d \big ) \to \mathcal B(\mathcal H)$$ such that  $\mathfrak D \big ( \big [\!\!\big[ \boldsymbol T^*, \boldsymbol T \big]\!\!\big] \big )$ is in the trace class. To achieve this, we introduce the determinant operator $\text{dEt}\,\big (\big[\!\!\big[{\boldsymbol M^{(\lambda)}}^* , \boldsymbol M^{(\lambda)}\big]\!\!\big]\big )$. For the weighted Bergman spaces $\mathcal H^{(\lambda)}(\mathbb B_d)$, we verify that the operator $\text{dEt}\,\big (\big[\!\!\big[{\boldsymbol T}^*,  \boldsymbol T\big]\!\!\big]\big )$ is in trace class. We also identify a large class $BS_{m,\vartheta}(\Omega)$ of commuting $d$ - tuple of operators $\boldsymbol T$ such that $\text{dEt}\,\big (\big[\!\!\big[\boldsymbol T^* , \boldsymbol T\big]\!\!\big]\big )$  is in trace class whenever $\boldsymbol T\in BS_{m,\vartheta}(\Omega)$.

\begin{defn}
Let $B_{i j}:\mathcal H \to \mathcal H$ be a bounded linear operator on the complex separable Hilbert space $\mathcal H$, $1\leq i,j \leq d$ and $\boldsymbol{B}:=\big (\! \big (B_{i j}\big )\! \big )$ be the  bounded linear operator from the Hilbert space $\mathcal H\otimes \ell_2(d)$ to itself. The \textit{determinant} $\text{dEt}\,(\boldsymbol B):\mathcal H\to \mathcal H$ is defined to be the operator given by the symmetrized version of  Laplace's expansion formula: 
$$\text{dEt}\,(\boldsymbol{B}):=\sum_{\sigma, \tau \in \mathfrak S_d} \text{Sgn}(\sigma)B_{\tau(1),\sigma(\tau(1))}B_{\tau(2),\sigma(\tau(2))} \ldots B_{\tau(d)\sigma(\tau(d))}.$$
For any commuting $d$- tuple $\boldsymbol{T}=(T_1,\ldots ,T_d)$, the determinant of the $d\times d$ block operator $\big[\!\! \big [\boldsymbol{T}^*, \boldsymbol{T}\big ]\!\!\big ]$ is  obtained by setting $B_{i j} = [T_j^*,T_i] $. \end{defn}

\subsection{Examples} In this subsection, we compute explicitly the determinant operator in several examples.  


Let $(M_{z_1}, M_{z_2})$ be the pair of multiplication operators on the Hardy space $H^2(\mathbb D^2)$ determined by the orthonormal basis $\{z_1^mz_2^n:m,n \geq 0\}$. The space $H^2(\mathbb D^2)$ is isometrically isomorphic to  $H^2(\mathbb{D})\otimes H^2(\mathbb{D})$ via the map $L : z_1^mz_2^n \mapsto z_1^m \otimes z_2^n$. Extend the map $L$ by linearity and note that it is well-defined and isometric. It is evidently surjective, hence unitary. The unitary $L$ intertwines the pair of operators 
$ (M_{z_1}, M_{z_2})$ with the pair $(M\otimes I, I \otimes M)$, where $M$ is the operator of multiplication by the coordinate function  on $H^2(\mathbb{D})$. We will let $\boldsymbol M$ denote either of these two pairs without causing any ambiguity since the meaning would be clear from the context.


The usual Hardy space $H^2(\mathbb D^2)$ is a module over the polynomial ring $\mathbb C[z_1,z_2]$ equipped with the module multiplication $m_p$ given by the point-wise multiplication, namely, 
$m_p(f) = p f$, $f\in H^2(\mathbb D^2)$, $p\in \mathbb C[z_1,z_2]$. 
Obviously, there are several other possibilities for the module multiplication. In this subsection, we consider a different module multiplication defined by the 
commuting pair of operators $\boldsymbol T=(T_1, T_2)$:
\[T_1=M \otimes I+I \otimes M\,\, \mbox{\rm and}\,\, T_2=M \otimes M\]
acting on the Hardy space $H^2(\mathbb{D}^2)$. The Hardy space equipped with the module multiplication: $m_p(f) = p(\boldsymbol T) f$, $f\in H^2(\mathbb D^2),\,p\in \mathbb C([z_1,z_2])$, is the Hardy module on the symmetrized bidisc. We have 
\begin{align*}
\big [\!\! \big [\boldsymbol{T}^*, \boldsymbol{T}\big ]\!\!\big ]=\begin{pmatrix}
P\otimes I & I\otimes P\\
P\otimes M & M\otimes P
\end{pmatrix}\begin{pmatrix}
P\otimes I & I\otimes P\\
P\otimes M & M\otimes P
\end{pmatrix}^*+\begin{pmatrix}
0 & 0\\
0 & P\otimes P
\end{pmatrix}\geq0.
\end{align*}
 It follows that $\boldsymbol{T}$ is hyponormal. A simple computation gives  $$\text{dEt}\,\big(\big[\!\! \big [\boldsymbol{T}^*, \boldsymbol{T}\big ]\!\!\big ]\big)=2  P\otimes P-P M^* \otimes M P- M P \otimes P M^*. $$
We note that the vector $1\otimes 1$ is an eigenvector of $\text{dEt}\,\big(\big[\!\! \big [\boldsymbol{T}^*, \boldsymbol{T}\big ]\!\!\big ]\big)$ with eigenvalue $2$ while the vector $z_1\otimes 1+1\otimes z_2$ is an eigenvector of $\text{dEt}\,\big(\big[\!\! \big [\boldsymbol{T}^*, \boldsymbol{T}\big ]\!\!\big ]\big)$ with eigenvalue $-1$. Therefore the operator $\text{dEt}\,\big(\big[\!\! \big [\boldsymbol{T}^*, \boldsymbol{T}\big ]\!\!\big ]\big)$ acting on $H^2(\mathbb{D}^2)$, is neither negative nor positive definite. However, we show that the restriction of 
$\text{dEt}\,\big(\big[\!\! \big [\boldsymbol{T}^*, \boldsymbol{T}\big ]\!\!\big ]\big)$ to the subspace $H_{\text{anti}}^2(\mathbb{D}^2)$ consisting of those functions in $H^2(\mathbb D^2)$ that are anti-symmetric, then it is nonnegative definite. 

Note that $H_{\text{anti}}^2(\mathbb{D}^2)$ is a reducing subspace for the pair $(T_1,T_2)$ and therefore the restriction of 
$\text{dEt}\,\big(\big[\!\! \big [\boldsymbol{T}^*, \boldsymbol{T}\big ]\!\!\big ]\big)$ to the subspace $H_{\text{anti}}^2(\mathbb{D}^2)$ equals $\text{dEt}\,\big(\big[\!\! \big [{\boldsymbol T}_{|\rm res}^*, \boldsymbol{T}_{|\rm res}\big ]\!\!\big ]\big)$, where $\boldsymbol T_{|\rm res}$ denotes the restriction of $\boldsymbol T$ to the subspace $H_{\text{anti}}^2(\mathbb{D}^2)$.

Let $[\!\![2]\!\!]$ be the set of all pairs $\boldsymbol{n}=(n_1, n_2)$ such that $n_1>n_2\geq 0,$ $n_1,n_2\in \mathbb N_0$.  Define
\[e_{\boldsymbol{n}}(\boldsymbol z):=\frac{z_1^{n_1}\otimes z_2^{n_2}-z_1^{n_2}\otimes z_2^{n_1}}{\sqrt{2}}.\]
Then $\{e_{\boldsymbol{n}}(\boldsymbol z):\boldsymbol{n}\in [\!\![2]\!\!]\}$ is an orthonormal basis for the subspace $H_{\text{anti}}^2(\mathbb{D}^2)$, see \cite{GSZ}. It is also shown  that $H_{\text{anti}}^2(\mathbb{D}^2)$ is module isomorphic to the Hardy module $H^2(G_2)$ on the symmetrized bidisc: $G_2:=\{(z_1+z_2, z_1 z_2): |z_1|, |z_2| < 1\}$. In other words, the multiplication by $p(T_1,T_2)$ on the Hardy space $H^2(\mathbb D^2)$ is unitarily equivalent to the multiplication by the pair of the coordinate functions on the Hardy space $H^2(G_2)$ of the symmetrized bidisc $G_2$. A direct and easy computation, using the orthonormal basis $\{e_{\boldsymbol n}\}$, shows that the   operator $\text{dEt}\,\big(\big[\!\! \big [\boldsymbol{T}^*, \boldsymbol{T}\big ]\!\!\big ]\big)$ is nonnegative definite and is in trace class:   
$$\big \langle \text{dEt}\,\big(\big[\!\! \big [\boldsymbol{T}^*, \boldsymbol{T}\big ]\!\!\big ]\big)e_{\boldsymbol n}, e_{\boldsymbol n}\big \rangle =  \begin{cases} 1 &  \mbox{\rm ~if~} \boldsymbol n =(1,0) \\ 0 & \mbox{\rm ~otherwise~}. \end{cases} $$

\subsection{Generalized Commutator} \label{subsec3.3}
Given any $d$- tuple of operators $\boldsymbol A$, not necessarily commuting, it is not clear what represents the degree of commutativity among these operators. For two operators $T_1, T_2$, the answer is clear, namely, one might demand that the commutator $[T_1,T_2]$ is small in an appropriate sense: some possibilities are zero, finite rank, trace class, compact. A particularly interesting special case occurs when $T
_2 = T_1^*$. The notion of a generalized commutator of a set of bounded operators $T_1, \ldots , T_{n}$ was introduced by Helton and Howe and may be thought of as a measure of the lack of commutativity among these operators. Before giving this definition, we recall a theorem of Amitsur–Levitzki \cite{AL}.
Let $S_h$ be the standard polynomial 
\begin{equation}\label{spol}
S_{h}(x_1, . . . , x_{h}) := \sum_{\sigma\in \mathfrak  S_{h}}\text{Sgn}(\sigma)x_{\sigma(1)}\cdots  x_{\sigma(h)}
\end{equation}
in non-commuting variables $x_1, \ldots ,x_h$. 
For any set of $2n$  element $A_1, \ldots , A_{2n}$ in the algebra of $n\times n$ matrices over a commutative ring, the Amitsur-Levitzki theorem asserts that $S_{2n}(A_1, \ldots , A_{2n})=0$,  see \cite{P}. The Helton and Howe proposed the following notion of a generalized commutator that has proved to be quite useful, (see \cite[Section A, p. 272]{HeltonHowe}).
 \begin{defn}[Helton-Howe]\label{Hel-Ho}
Let $\boldsymbol{A}= (A_1, \ldots , A_n)$ be a $n$- tuple of bounded operators. The generalized commutator
$\text{GC}\,(\boldsymbol A)$ is defined to be $S_n(A_1, \ldots, A_n)$. 
\end{defn}

 We adapt the definition of Helton and Howe slightly to the case of a commuting $d$- tuple of operators $\boldsymbol T$ as follows.  Let $\boldsymbol{T}=(T_1,\ldots ,T_d)$ be a $d$- tuple of operators. Let $A_1=T_1^*,A_2=T_1,\ldots , A_{2d-1}=T_d^*, A_{2d}=T_d$. The generalized commutator $\text{GC}(\boldsymbol{T}^*, \boldsymbol{T})$ is defined to be the sum  
\begin{equation}\label{GC(T*,T)}
\text{GC}(\boldsymbol{T}^*,\boldsymbol{T}):= S_{2d}(A_1, A_2, \ldots , A_{2d-1},A_{2d}).
\end{equation}

Recall, Definition 2.2 of \cite{S}, that a $n$ - normal operator $T$ is a bounded linear operator on a Hilbert space $\mathcal H$ which is unitarily equivalent to an operator of the form
$$A[n]:= \big(\!\big( A_{ij}\big)\!\big):\widetilde{ \mathcal H }\otimes \mathbb C^n \to \widetilde{\mathcal H} \otimes \mathbb C^n,$$
for some family $A_{ij}:\widetilde{\mathcal H} \to \widetilde{\mathcal H}$, $1\leq i,j \leq n$, of commuting normal operators.  
In an analogous manner, we say that a commuting $d$ - tuple $\boldsymbol T:= (T_1, T_2, \ldots , T_d)$ is $n$- normal if $T_k$ is unitarily equivalent to an operator of the form
$$A_k[n]:= \big(\!\big( A^{(k)}_{ij}\big)\!\big): \widetilde{\mathcal H} \otimes \mathbb C^n \to \widetilde{\mathcal H} \otimes \mathbb C^n,$$
where the operator $A_k[n]$ is $n$ - normal and 
$A_k[n]$, $1\leq k\leq d,$ is commuting. 

If $\boldsymbol T$ is $d$- normal, then applying the Amitsur-Levitzki theorem to the $d$ - tuple $\boldsymbol A[d]:=\big (A_1[d],\ldots, A_d[d]\big )$, we conclude that $\text{GC}(\boldsymbol A[d]^*,\boldsymbol A[d])=0$. Thus 
\begin{equation}\label{A-L}
\text{GC}(\boldsymbol T^*,\boldsymbol T)= U[d]^*\big(\text{GC}(\boldsymbol A[d]^*,\boldsymbol A[d])\big ) U[d]=0,\end{equation}
where the unitary $U[d]: \mathcal H \to \widetilde{\mathcal H}\otimes \mathbb C^d,$ intertwining $T_k$ and $A_k[d],$ is independent of $k=1,\ldots , d$.


We now show that the  $\text{dEt}\,\big(\big[\!\! \big [\boldsymbol{T}^*, \boldsymbol{T}\big ]\!\!\big ]\big)$ and $\text{GC}\,(\boldsymbol T^*, \boldsymbol T)$ coincide for any commuting tuple $\boldsymbol T$. We emphasize that the equality need not hold unless $\boldsymbol T$ is a $d$ - tuple of commuting operators. In this case, working with the $\text{dEt}\,\big(\big[\!\! \big [\boldsymbol{T}^*, \boldsymbol{T}\big ]\!\!\big ]\big)$ has some advantages over $\text{GC}\,(\boldsymbol T^*, \boldsymbol T)$ since a number of terms in $\text{GC}\,(\boldsymbol T^*, \boldsymbol T)$ cancel, in case $\boldsymbol T$ is $d$ - tuple of commuting operators, and it equals the less formidable expression for $\text{dEt}\,\big(\big[\!\! \big [\boldsymbol{T}^*, \boldsymbol{T}\big ]\!\!\big ]\big)$. 
\begin{thm} \label{HH} For any d-tuple $\boldsymbol{T}$ of commuting operators, the determinant 
$$\mbox{\rm dEt}\,\big(\big[\!\! \big [\boldsymbol{T}^*, \boldsymbol{T}\big ]\!\!\big ]\big)= \text{GC}\,(\boldsymbol T^*, \boldsymbol T).$$ 
\end{thm}
\begin{proof} By definition, we have  
	\begin{align*}
	\text{dEt}\,\big(\big[\!\! \big [\boldsymbol{T}^*, \boldsymbol{T}\big ]\!\!\big ]\big) &=\sum_{\tau ,\sigma \in \mathfrak  S_d}\text{Sgn}(\sigma)\prod_{i=1}^{d}B_{\tau(i)\sigma(\tau(i))}\\
	&=\sum_{\tau ,\eta \in \mathfrak  S_d}\text{Sgn}(\tau)\text{Sgn}(\eta)\prod_{i=1}^{d}B_{\tau(i)\eta(i)}, 
	\end{align*} 
	where $B_{i j} = [T_j^*,T_i]$ and $\eta=\sigma\tau$. 
	
	Fix a commuting tuple of operators $\boldsymbol T$. Suppose one of the terms in  $\text{GC}(\boldsymbol{T^*, \boldsymbol T})$ has a string of the form $P T_i T_j Q$, where $P$ and $Q$ are products of operators taken from the remaining set of $(2d-2)$ operators: $(T_1^*, \ldots, T_d^*, T_1, \ldots ,\hat{T}_i, \ldots, \hat{T}_j, \ldots ,T_d)$. (Here $i,j$ are from $\{1,2, \ldots ,d\}$ and $\hat{T}$ means that it is not included in the set.)
	Then there must be a second term in $\text{GC}(\boldsymbol{T^*, \boldsymbol T})$  of the form  $P T_j T_i Q$ with the opposite sign. However these have to cancel since $T_i T_j = T_j T_i$.  A similar argument applies to strings of the form $R T_i^*T_j^* S$. Thus the only terms that survive are those in which a $T_i$ must be followed by a $T_j^*$ and a $T_j^*$ must be followed by a $T_i$.
	
	There are two sets of terms in $\text{GC}(\boldsymbol{T}^*, \boldsymbol T)$, one set which begins with a $T^*$ and another set which begins with a $T$. Indeed, we have
	\begin{align}\label{cherian}
	\text{GC}(\boldsymbol{T}^*, \boldsymbol{T})= &\sum_{\tau,\eta\in \mathfrak  S_d}\text{Sgn}(\tau)\text{Sgn}(\eta)T^*_{\eta(1)}T_{\tau(1)}T^*_{\eta(2)}\ldots T^*_{\eta(d)}T_{\tau(d)}\nonumber\\
	&+ (-1)^d\sum_{\tau,\eta\in \mathfrak  S_d}\text{Sgn}(\tau)\text{Sgn}(\eta)T_{\tau(1)}T^*_{\eta(1)}T_{\tau(2)}\ldots T_{\tau(d)}T^*_{\eta(d)}.  \end{align}
	The terms starting with a $T^*$, which is the first sum in $\text{GC}(\boldsymbol{T}^*, \boldsymbol T)$ simplifies:
	\begin{align*}& \sum_{\tau,\eta\in \mathfrak S_d}\text{Sgn}(\tau)\text{Sgn}(\eta) { T^*_{\eta(1)}T_{\tau(1)}T^*_{\eta(2)}\ldots T^*_{\eta(d)}T_{\tau(d)} } \\&  =\sum_{\tau,\eta\in \mathfrak S_d}\text{Sgn}(\tau)\text{Sgn}(\eta){ \big [ T^*_{\eta(1)}, T_{\tau(1)}\big ] T^*_{\eta(2)}\ldots T^*_{\eta(d)}T_{\tau(d)} } \\& + \sum_{\tau,\eta\in \mathfrak S_d}\text{Sgn}(\tau)\text{Sgn}(\eta) { T_{\tau(1)}T^*_{\eta(1)} T^*_{\eta(2)}\ldots T^*_{\eta(d)}T_{\tau(d)} }.
	\end{align*}
	If $d\geq1$ the second sum on the right is zero since there are two terms containing the string $T^*_{\eta(1)} T^*_{\eta(2)}$ with opposite signs.  Repeating this process (note that the vanishing argument does not apply at the last stage), we get
	\begin{align}\label{firstsum}
	\sum_{\tau,\eta\in \mathfrak S_d} & \text{Sgn}(\tau)  \text{Sgn}(\eta) { T^*_{\eta(1)}T_{\tau(1)}T^*_{\eta(2)}\ldots T^*_{\eta(d)}T_{\tau(d)} } \nonumber\\&=\sum_{\tau,\eta\in \mathfrak S_d}\text{Sgn}(\tau)\text{Sgn}(\eta) \big[ T^*_{\eta(1)}, T_{\tau(1)}\big] \big[ T^*_{\eta(2)}, T_{\tau(2)}\big]\ldots \big[ T^*_{\eta(d)}, T_{\tau(d)}\big]\nonumber \\
	+\sum_{\tau,\eta\in \mathfrak S_d} & \text{Sgn}(\tau)\text{Sgn}(\eta) \big[ T^*_{\eta(1)}, T_{\tau(1)}\big] \big[ T^*_{\eta(2)}, T_{\tau(2)}\big]\ldots \big[ T^*_{\eta(d-1)}, T_{\tau(d-1)}\big]T_{\tau(d)}T^*_{\eta(d)}.  \end{align}
			The terms starting with a $T$, which is the second sum in $\text{GC}(\boldsymbol{T}^*, \boldsymbol T)$, using the equality 
	 	\begin{align*}
	T_{\tau(i)}T^*_{\eta(i)}=-[T^*_{\eta(i)},T_{\tau(i)}]+T^*_{\eta(i)}T_{\tau(i)},
	\end{align*}
	simplifies:
	\begin{align*}
	& \sum_{\tau,\eta\in \mathfrak S_d}\text{Sgn}(\tau)\text{Sgn}(\eta) { T_{\tau(1)}T^*_{\eta(1)}T_{\tau(2)}\ldots T_{\tau(d)}T^*_{\eta(d)} } \\&  =-\sum_{\tau,\eta\in \mathfrak S_d}\text{Sgn}(\tau)\text{Sgn}(\eta) \big[ T^*_{\eta(1)}, T_{\tau(1)}\big] T_{\tau(2)}\ldots T_{\tau(d)}T^*_{\eta(d)}  \\& + \sum_{\tau,\eta\in \mathfrak S_d}\text{Sgn}(\tau)\text{Sgn}(\eta)  T^*_{\eta(1)}T_{\tau(1)} T_{\tau(2)}\ldots T_{\tau(d)}T^*_{\eta(d)} .
	\end{align*}
	If $d\geq1$ the second sum on the right is zero since there is a string with $T_{\tau(1)} T_{\tau(2)}$. Repeating this process $d-1$ times, we get
	\begin{align}\label{secondsum}
	(-1)^d &\sum_{\tau,\eta\in \mathfrak S_d}\text{Sgn}(\tau)\text{Sgn}(\eta) { T_{\tau(1)}T^*_{\eta(1)}T_{\tau(2)}\ldots T_{\tau(d)}T^*_{\eta(d)} } \nonumber \\ = (-1)^d (-1)^{d-1} &   \sum_{\tau,\eta\in \mathfrak S_d}\text{Sgn}(\tau)\text{Sgn}(\eta) \big[ T^*_{\eta(1)}, T_{\tau(1)}\big] \big[ T^*_{\eta(2)}, T_{\tau(2)}\big]\ldots \big[ T^*_{\eta(d-1)}, T_{\tau(d-1)}\big]T_{\tau(d)}T^*_{\eta(d)}.
	\end{align}
	Adding the two sums on the right hand side of the equation  \eqref{firstsum} and the one on the right hand side of the equation \eqref{secondsum}, we get
	\begin{align*}
	\text{GC}(\boldsymbol{T}^*, \boldsymbol{T})=&\sum_{\tau,\eta\in \mathfrak S_d}\text{Sgn}(\tau)\text{Sgn}(\eta) [ T^*_{\eta(1)}, T_{\tau(1)}][ T^*_{\eta(2)}, T_{\tau(2)}]\ldots [ T^*_{\eta(d)}, T_{\tau(d)}]\\&=\sum_{\tau,\eta\in \mathfrak S_d}\text{Sgn}(\tau)\text{Sgn}(\eta){B_{\tau(1)\eta(1)}B_{\tau(2)\eta(2)}\ldots B_{\tau(d)\eta(d)}}\\&=\text{dEt}\,\big(\big[\!\! \big [\boldsymbol{T}^*, \boldsymbol{T}\big ]\!\!\big ]\big)
	\end{align*}
completing the verification that $GC(\boldsymbol{T}^*, \boldsymbol{T}) = \text{dEt}\big(\big[\!\! \big [\boldsymbol{T}^*, \boldsymbol{T}\big ]\!\!\big ]\big)$ for a commuting tuple $\boldsymbol{T}$.	
\end{proof}

\begin{cor} \label{cor:3.5}
Let $\boldsymbol N=(N_1, N_2)$ be a commuting pair of normal operators and $\boldsymbol A=(A_1, A_2)$ be a pair of $2\times 2$ scalar matrices. If $\boldsymbol N \otimes \boldsymbol A:=(N_1\otimes A_1, N_2\otimes A_2)$ is a $2$ - normal pair, then $\text{dEt}\,\big(\big[\!\!\big[(\boldsymbol N \otimes \boldsymbol A)^*, \boldsymbol N \otimes \boldsymbol A \big]\!\!\big]\big)=0.$
\end{cor}

\begin{proof} We have $\text{dEt}\,\big(\big[\!\!\big[(\boldsymbol N \otimes \boldsymbol A)^*, \boldsymbol N \otimes \boldsymbol A \big]\!\!\big]\big)= N_1 N_2 N_1^* N_2^* \otimes \text{dEt}\,\big(\big[\!\!\big[\boldsymbol A^*, \boldsymbol A \big]\!\!\big]\big)$. As $\boldsymbol N \otimes \boldsymbol A$ is $2$ - normal, we also have  $[N_1\otimes A_1, N_2\otimes A_2]=0$ implying  
either $N_1N_2=0$ or $[A_1, A_2]=0.$
If $N_1 N_2=0$, then there is nothing to prove. If $[A_1, A_2]=0$, then $\text{dEt}\,\big(\big[\!\!\big[\boldsymbol A^*, \boldsymbol A \big]\!\!\big]\big)= \text{GC}\,(\boldsymbol A^*, \boldsymbol A)$ by Thorem \ref{HH} and the Amitsur-Levitzki theorem shows that  $\text{GC}\,(\boldsymbol A^*, \boldsymbol A)=0.$ 
\end{proof}
This corollary is easily strengthened using  the full force of the Amitsur-Levitzki theorem together with the equality $GC(\boldsymbol T^*, \boldsymbol T) = \text{dEt}\,\big(\big[\!\!\big[(\boldsymbol T^*, \boldsymbol T \big]\!\!\big]\big)$. 
\begin{cor}\label{cor:3.6}
If a commuting $d$- tuple $\boldsymbol T$ is $d$- normal, then $\text{dEt}\,\big(\big[\!\! \big [\boldsymbol{T}^*, \boldsymbol{T}\big ]\!\!\big ]\big ) = 0.$
\end{cor}
\begin{proof}
    First, note that  $\text{GC}(\boldsymbol T^*,\boldsymbol T) = 0$ for any $d$- tuple of $d$- normal operators, see Equation \eqref{A-L}. Combining this with Theorem \ref{HH} asserting the equality of $\text{GC}(\boldsymbol T^*, \boldsymbol T)$ and  $\text{dEt}\,\big(\big[\!\! \big [\boldsymbol{T}^*, \boldsymbol{T}\big ]\!\!\big ]\big )$, we have the desired conclusion.
\end{proof}

For any pair $N_1,N_2$ of commuting normal operators acting on a Hilbert space $\mathcal H$, consider the operators $T_i:\mathcal H\otimes \mathbb C^2 \to \mathcal H \otimes \mathbb C^2$ of the form: 
$T_i=\left ( \begin{smallmatrix} N_i& \lambda_i I\\ 0& N_i\end{smallmatrix}\right )$, $i=1,2$. The following properties of the pair $\boldsymbol T:=(T_1,T_2)$ are easily verified. 
\begin{itemize}
\item The pair $\boldsymbol T$ is commuting and $2$ - normal.
    \item The pair $\boldsymbol T$ is not hyponormal, in general. For instance, take $\lambda_1=1=\lambda_2$. 
\item  $\text{dEt}\,\big(\big[\!\! \big [\boldsymbol{T}^*, \boldsymbol{T}\big ]\!\!\big ]\big )=0$.  
\end{itemize}
This shows that even if $\text{dEt}\,\big(\big[\!\! \big [\boldsymbol{T}^*, \boldsymbol{T}\big ]\!\!\big ]\big )=0$, the pair of commuting operators $(T_1,T_2)$ need not be hyponormal. 
\section{Trace estimate of the determinant operator in the class $BS_{m,\vartheta}(\Omega)$} Let $\Omega \subset \mathbb{C}^d$ be a bounded domain.  In this section, we define a class  $BS_{m,\vartheta}(\Omega)$ consisting of commuting $d$- tuple of operators acting on  a Hilbert space $\mathcal{H}$. 
Our main result  is an estimate for the trace of $\text{dEt}\,\big(\big[\!\! \big [\boldsymbol{T}^*, \boldsymbol{T}\big ]\!\!\big ]\big )$ for $d$- tuples of operators $\boldsymbol T$ in $BS_{m,\vartheta}(\Omega)$. The proof is modeled  after the first part of Voiculescu's proof of the Berger-Shaw theorem in \cite{Voiculescu}.
The finite dimensional subspaces of $\mathcal H$ defined below play a significant role in our computation of the trace. 
 
\begin{defn} \label{P_N} For a $m$- cyclic $d$- tuple $\boldsymbol{T}$, let 
\[ \mathcal{H}_N:=\bigvee \bigg\{T_{1}^{i_1}T_{2}^{i_2}\ldots T_{d}^{i_d}v |~ v\in \boldsymbol{\xi}[m] ~\text{and}~ 0\leq i_{1}+i_{2}+\ldots i_{d} \leq N \bigg \} \] and $P_{N}$ be the orthogonal projection onto $\mathcal{H}_N$. 
\end{defn}

We list below some of the basic properties of the projection $P_N$  that will be used in the proof of the main theorem.

\begin{lem}\label{rank} For a $m$-cyclic $d$-tuple of operators $\boldsymbol T$, we have $P_N$ increasing strongly to $I$ and  $\text{rank}\,(P_N ^{\perp}T_j P_N)\leq m\binom{N+d-1}{d-1}$,\,\, $j =1,\ldots, d$.
\end{lem}
\begin{proof}
Evidently, $\mathcal H_N \subseteq \mathcal H_{N+1}$ and hence 
the projections $\{P_N\}$ are increasing. By hypothesis, $\boldsymbol T$ is $m$-cyclic, therefore by definition, the linear span of $\{\mathcal H_N:N\in \mathbb N_0\}$ is dense in $\mathcal H $.  
The number of vectors from $\mathcal H_N$ that are pushed out of it by the operator $T_j$ provides a reasonable upper bound on the rank of the operator $P_N^{\perp}T_j P_N$. Such vectors can only be a subset of the subspace $ \mathcal H_N\ominus \mathcal H_{N-1}$. 
Clearly, the dimension of this subspace is the same as the dimension of the space of homogeneous polynomials of degree $N$ in $d$-variables tensored with $\mathbb C^m$, which is $m \binom{N+d-1}{d-1}.$ Therefore,  $\text{rank}\,(P_N ^{\perp}T_j P_N)\leq m\binom{N+d-1}{d-1}.$
\end{proof}
We recall a well-known inequality between the trace norm and the operator norm of a finite rank bounded operator. 
\begin{lem}
If $F\in \mathcal{B}(\mathcal{H})$ is of finite rank, then  $\|F\|_1 \leq (\text{rank}\, F) \|F\|$.
\end{lem}

\begin{proof}
Since $F$ is a finite rank operator, choosing an arbitrary but fixed   $\{\varphi_1, \ldots ,\varphi_n\}$ orthonormal basis for the range of $F$, we have  
\begin{equation}\label{FRTN}
F x = \sum_{k=1}^n \langle x, v_k\rangle\varphi_k,\,\, x\in \mathcal H \end{equation}
for some set of $n$ vectors $\{v_1,\ldots ,v_n\}$ in $\mathcal H$. 
This representation of the operator $F$ shows that  $F^* y = \big ( \sum_{k=1}^n \langle \varphi_k , y \rangle v_k\big ),$ $y\in \mathcal H$. Consequenly, 
\begin{align*}
 F^*F x &= \sum_{j=1}^n \langle x, v_j \rangle v_j,\,\, x\in \mathcal H.
\end{align*}
Therefore, setting $V$ to be the linear span of the vectors $\{v_1, \ldots, v_n\}$, we conclude that $F^*F = (F^* F)_{|V} \oplus 0$. Since $V$ is finite dimensional,  it follows that 
$$\|F\|_1 = \big \|(F^* F)^{1/2}\big \|_1 \leq (\text{rank}\, F) \|F\|$$
completing the proof of the lemma.
\end{proof}

We next define a class $BS_{m,\vartheta}(\Omega)$ of  $d$ - tuples of commuting operators. The rest of this section is devoted to showing that if $\boldsymbol T$ is in  $BS_{m,\vartheta}(\Omega)$, then $\text{trace} \big (\text{dEt}\,\big (\big[\!\! \big [\boldsymbol{T}^*, \boldsymbol{T}\big ]\!\!\big ]\big) \big )$ is finite. Moreover, if $\boldsymbol T$ is in $BS_{m,\vartheta}(\Omega)$, then 
$$\text{trace} \big (\text{dEt}\,\big (\big[\!\! \big [\boldsymbol{T}^*, \boldsymbol{T}\big ]\!\!\big ]\big) \big )\leq m\, \vartheta \,d!\prod_{i=1}^{d}\|T_i\|^2.$$


\begin{defn} \label{class} Fix a bounded domain $\Omega\subset \mathbb C^d$ such that $\overbar{\Omega}$ is polynomially convex. 
A  $m$-cyclic commuting $d$ - tuple of operators with  $\sigma(\boldsymbol{T})=\overline{\Omega}$ is said to be in the class $BS_{m, \vartheta}(\Omega)$, if 
\begin{itemize}
		\item[\rm (i)] $P_N T_j P_N^\perp = 0$, $j=1,\ldots , d$.  
		\item[\rm (ii)] $\text{dEt}\,\big (\big[\!\! \big [\boldsymbol{T}^*, \boldsymbol{T}\big ]\!\!\big ]\big)$ is non-negative definite.
		\item[\rm (iii)]For a fixed but arbitrary $\tau$ in the permutation group $\mathfrak S_{d}$ of $d$ symbols, there exists $\vartheta \in \mathbb N$, independent of $N$, such that $$\big \|P_N \big ( \sum_{\eta\in \mathfrak S_d}\text{Sgn}(\eta) T^*_{\eta(1)}T_{\tau(1)}T^*_{\eta(2)}\ldots T^*_{\eta(d)}\big ) P_N^{\perp} T_{\tau(d)}P_N\big\|\leq \vartheta\, {\binom{N+d-1}{d-1}}^{-1}\prod_{i=1}^d \big \|T_i\big \|^2 .$$
	\end{itemize}
\end{defn}
\begin{rem} 
\begin{itemize} 
\label{general tau}
\item[]
 \item[\rmfamily(a)] For a single operator T on a Hilbert space $\mathcal{H}$, condition {\rm (iii)} of Definition \ref{class} reduces to  
$$\|P_N T^* P_N^{\perp} T P_N\|\leq  \vartheta \|T\|^2,$$
which is true with $\vartheta = 1$. It follows that a $m$-cyclic hyponormal operator $T$ with $\sigma(T)=\overline{\Omega}$ is in the class $BS_{m, 1}(\Omega)$, if $P_N T P_N^\perp = 0$. 
\item[\rmfamily(b)]
If for each $\tau \in \mathfrak S_d$, there exists a unitary operator $U_{\tau}$ on the Hilbert space such that 
$$U_{\tau}T_{\tau(i)}U_{\tau}^*=T_i,\,\, 1 \leq i \leq d,$$ 
then it is enough to check condition (iii) for identity permutation, that is,
$$\|P_N \big (\sum_{\eta \in \mathfrak S_{d}}\text{Sgn}(\eta)T_{\eta(1)}^*T_{1}T_{\eta(2)}^*\ldots T_{d-1} T_{\eta(d)}^*\big )P_N^{\perp} T_{d} P_N\|\leq \vartheta {\binom{N+d-1}{d-1}}^{-1}\prod_{i=1}^d \|T_i\|^2 $$
implies all the other inequalities, one for each $\tau$ of (iii) in Definition \ref{class}. 

To see this, pick $\tau_0\in \mathfrak S_d$ such that 
$\eta=\tau \cdot \tau_0$. With this choice of $\tau_0$, we have $\text{Sgn}\,\eta=\text{Sgn}\,\tau \text{Sgn}\,\tau_0$. It now follows that 
\begin{align*}
 P_N &\big (\sum_{\eta \in \mathfrak S_{d}}\text{Sgn}(\eta)T_{\eta(1)}^*T_{\tau(1)}T_{\eta(2)}^*\ldots T_{\tau(d-1)} T_{\eta(d)}^*\big )P_N^{\perp} T_{\tau(d)} P_N\\&=P_N U_{\tau}\big (\text{Sgn}(\tau)\sum_{\hat{\eta} \in \mathfrak S_{d}}\text{Sgn}(\hat{\eta})T_{\hat{\eta}(1)}^*T_{1}T_{\hat{\eta}(2)}^*\ldots T_{d-1} T_{\hat{\eta}(d)}^*\big )U_{\tau}^*P_N^{\perp} U_{\tau}T_{d}U_{\tau}^*P_N\\&=\text{Sgn}(\tau)U_{\tau}P_N \big (\sum_{\hat{\eta} \in \mathfrak S_{d}}\text{Sgn}(\hat{\eta})T_{\hat{\eta}(1)}^*T_{1}T_{\hat{\eta}(2)}^*\ldots T_{d-1} T_{\hat{\eta}(d)}^*\big )P_N^{\perp} T_{d}P_N U_{\tau}^*.
\end{align*}
\item[\rmfamily(c)] There exists unitary representation $U$ of the symmetric group $\mathfrak S_d$ and commuting $d$ - tuples of operators $\boldsymbol T$ in $BS_{1,1}(\mathbb B_d)$ with the property  $$U_{\tau}T_{\tau(i)}U_{\tau}^*=T_i,\,\, \tau \in \mathfrak S_n ,\,\, 1 \leq i \leq d.$$ 
\end{itemize}
\end{rem}
Explicit examples are given in the following section. We begin by proving a lemma, which is the main ingredient in the proof of the trace inequality of Theorem \ref{lem63}.


\begin{lem} \label{L1.6}
Assume that the $d$ -tuple $\boldsymbol{T}$ is $m$- cyclic and that $P_N T_j P_N^\perp = 0$, $1\leq j \leq d$. Then \begin{multline*}
\big |\text{trace}\,\big(P_N\text{dEt}\,\big(\big[\!\! \big [\boldsymbol{T}^*, \boldsymbol{T}\big ]\!\!\big ]\big)P_N\big)\big|\\ \leq m\,\tbinom{N+d-1}{d-1} \Big(\sum_{\tau\in \mathfrak S_d} \big \|\big(\sum_{\eta\in \mathfrak S_d}\text{Sgn}(\eta)P_N T^*_{\eta(1)}T_{\tau(1)}T^*_{\eta(2)}\ldots T^*_{\eta(d)}P_N^{\perp}T_{\tau(d)}P_N\big)\big \|\Big). \end{multline*}
\end{lem}

\begin{proof}
  For a d-tuple of commuting operators $\boldsymbol{T}$, by Proposition \ref{HH} and using Equation \eqref{cherian}, we infer that the determinant 
 \begin{align}\label{det}
  \text{dEt}\big(\big[\!\! \big [\boldsymbol{T}^*, \boldsymbol{T}\big ]\!\!\big ]\big)
     =&\sum_{\tau,\eta\in \mathfrak S_d}\text{Sgn}(\tau)\text{Sgn}(\eta)T^*_{\eta(1)}T_{\tau(1)}T^*_{\eta(2)}\ldots T^*_{\eta(d)}T_{\tau(d)}+\nonumber\\
     &\phantom{Paramita} (-1)^d\sum_{\tau,\eta\in \mathfrak S_d}\text{Sgn}(\tau)\text{Sgn}(\eta)T_{\tau(1)}T^*_{\eta(1)}T_{\tau(2)}\ldots T_{\tau(d)}T^*_{\eta(d)}\nonumber\\
      =&\sum_{\tau,\eta\in \mathfrak S_d}\text{Sgn}(\tau)\text{Sgn}(\eta)T^*_{\eta(1)}T_{\tau(1)}T^*_{\eta(2)}\ldots T^*_{\eta(d)}T_{\tau(d)}+\nonumber\\
      & \phantom{Paramita}-\sum_{\tau,\eta\in \mathfrak S_d}\text{Sgn}(\tau)\text{Sgn}(\eta)T_{\tau(d)}T^*_{\eta(1)}T_{\tau(1)}\ldots T_{\tau(d-1)}T^*_{\eta(d)}\\
      =&\sum_{\tau,\eta\in \mathfrak S_d}\text{Sgn}(\tau)\text{Sgn}(\eta)\big[T^*_{\eta(1)}T_{\tau(1)}T^*_{\eta(2)}\ldots T^*_{\eta(d)}, T_{\tau(d)}\big].\nonumber
 \end{align}
 The second equality in Equation \eqref{det} is obtained by replacing the permutation $\tau$ in the second sum by the permutation $\tau^\prime= \tau\circ (1,\ldots , d)$, where $(1,\ldots ,d)$ is the permutation taking $1\to 2$, \ldots, $(d-1) \to d$ and $d \to 1$. Since the sums are over all permutations, this does not change it. But the sign of $\tau^\prime$ differs from that of $\tau$ by $(-1)^{d-1}$. Therefore, it follows that 
\begin{align*}
P_N\, \text{dEt}\,\big(\big[\!\! \big [\boldsymbol{T}^*, \boldsymbol{T}\big ]\!\!\big ]\big) P_N =&\sum_{\tau,\eta\in S_d}\text{Sgn}(\tau)\text{Sgn}(\eta)P_N\big[T^*_{\eta(1)}T_{\tau(1)}T^*_{\eta(2)}\ldots T^*_{\eta(d)}, T_{\tau(d)}\big]P_N.
\end{align*}
Also, for $\tau \in \mathfrak S_d$, 
\begin{align*}
 P_N\big[T^*_{\eta(1)}T_{\tau(1)}T^*_{\eta(2)}\ldots T^*_{\eta(d)}, T_{\tau(d)}\big]P_N \phantom{Par}&\\ 
 = P_N\big (T^*_{\eta(1)}T_{\tau(1)}T^*_{\eta(2)}\ldots T^*_{\eta(d)}T_{\tau(d)}
 & - T_{\tau(d)}T^*_{\eta(1)}T_{\tau(1)}\ldots T_{\tau(d-1)}T_{\eta(d)}\big )P_N \\
 =P_N T^*_{\eta(1)}T_{\tau(1)}T^*_{\eta(2)}\ldots T^*_{\eta(d)}(P_N+P_N^{\perp})T_{\tau(d)} P_N
 &-P_NT_{\tau(d)}(P_N+P_N^{\perp})T^*_{\eta(1)}T_{\tau(1)}T^*_{\eta(2)}\ldots T^*_{\eta(d)}P_N\\
 = P_N T^*_{\eta(1)}T_{\tau(1)}T^*_{\eta(2)}\ldots T^*_{\eta(d)}P_N^{\perp}T_{\tau(d)}P_N &+ \big [P_NT^*_{\eta(1)}T_{\tau(1)}T^*_{\eta(2)}\ldots T^*_{\eta(d)}P_N, P_N T_{\tau(d)} P_N\big ] \\  -P_NT_{\tau(d)}P_N^{\perp}T^*_{\eta(1)}T_{\tau(1)}T^*_{\eta(2)}\ldots& T^*_{\eta(d)}P_N\\
 = P_N T^*_{\eta(1)}T_{\tau(1)}T^*_{\eta(2)}\ldots T^*_{\eta(d)}P_N^{\perp}T_{\tau(d)}P_N &+ \big [P_NT^*_{\eta(1)}T_{\tau(1)}T^*_{\eta(2)}\ldots T^*_{\eta(d)}P_N, P_N T_{\tau(d)} P_N\big ],
\end{align*}
where, the validity of the last equality follows from the assumption that  $P_N T_j P_N^{\perp}=0$, $j=1,2,\ldots ,d$.
If $A, B$ are any two operators with one in trace class and the other bounded, then $\text{trace}\,(AB)=\text{trace}\,(BA)$.
A bounded operator of finite rank is trace class. Therefore, the commutator of any two  bounded operators of finite rank must be $0$. Hence
$$\text{trace}\,\Big (\big [P_NT^*_{\eta(1)}T_{\tau(1)}T^*_{\eta(2)}\ldots T^*_{\eta(d)}P_N, P_N T_{\tau(d)} P_N\big ]\Big )=0.$$
Putting these together we obtain the equality below
\begin{multline*}
 \text{trace}\,(P_N\, \text{dEt}\,\big(\big[\!\! \big [\boldsymbol{T}^*, \boldsymbol{T}\big ]\!\!\big ]\big) P_N)\\=\sum_{\tau \in \mathfrak S_d}\text{Sgn}(\tau)\text{trace}\big(\sum_{\eta\in \mathfrak S_d}\text{Sgn}(\eta)P_N T^*_{\eta(1)}T_{\tau(1)}T^*_{\eta(2)}\ldots T^*_{\eta(d)}P_N^{\perp}T_{\tau(d)}P_N\big).  
\end{multline*}
We also have the inequalities:
\begin{align}\label{rankineq}
\text{rank}\,\Big(\big(\sum_{\eta\in \mathfrak S_d}\text{Sgn}(\eta)P_N T^*_{\eta(1)}T_{\tau(1)}T^*_{\eta(2)} T^*_{\eta(d)}& P_N^{\perp}\big)\big (P_N^{\perp}T_{\tau(d)}P_N\big)\Big)\nonumber\\&\leq \text{rank} (P_N^{\perp}T_{\tau(d)}P_N) \leq  m\tbinom{N+d-1}{d-1}.
\end{align}
This implies
\begin{align} \label{trace}
\big |\text{trace}\,(P_N\, \text{dEt}&\,\big(\big[\!\! \big [\boldsymbol{T}^*, \boldsymbol{T}\big ]\!\!\big ]\big)  P_N)\big|\nonumber \\& \leq  \sum_{\tau \in \mathfrak S_d} \big |\text{trace}\,\big(\sum_{\eta\in \mathfrak S_d}\text{Sgn}(\eta)P_N T^*_{\eta(1)}T_{\tau(1)}T^*_{\eta(2)}\ldots T^*_{\eta(d)}P_N^{\perp}T_{\tau(d)}P_N\big) \big |\nonumber \\
& \leq \sum_{\tau \in \mathfrak S_d}\Big\{ \big\|\big(\sum_{\eta\in \mathfrak S_d}\text{Sgn}(\eta)P_N T^*_{\eta(1)}T_{\tau(1)}T^*_{\eta(2)}\ldots T^*_{\eta(d)}P_N^{\perp}T_{\tau(d)}P_N\big)\big \|\nonumber \\ &\hspace{7em}\text{rank}\,\big(\sum_{\eta\in \mathfrak S_d}\text{Sgn}(\eta)P_N T^*_{\eta(1)}T_{\tau(1)}T^*_{\eta(2)}\ldots T^*_{\eta(d)}P_N^{\perp}T_{\tau(d)}P_N\big)\Big\} \nonumber\\ &\leq \sum_{\tau\in \mathfrak S_d} \big \|\big(\sum_{\eta\in \mathfrak S_d}\text{Sgn}(\eta)P_N T^*_{\eta(1)}T_{\tau(1)}T^*_{\eta(2)}\ldots T^*_{\eta(d)}P_N^{\perp}T_{\tau(d)}P_N\big)\big \| \text{rank}\,\big(P_N ^{\perp}T_{\tau(d)} P_N\big)\nonumber \\
& \leq m\tbinom{N+d-1}{d-1}\Big(\sum_{\tau\in \mathfrak S_d}\big \|\big(\sum_{\eta\in \mathfrak S_d}\text{Sgn}(\eta)P_N T^*_{\eta(1)}T_{\tau(1)}T^*_{\eta(2)}\ldots T^*_{\eta(d)}P_N^{\perp}T_{\tau(d)}P_N\big)\big \|\Big)\nonumber.
\end{align}
The two penultimate inequalities follow from the inequality \ref{rankineq}.
\end{proof}

The following Theorem shows that the operator $\text{dEt}\big(\big[\!\! \big [\boldsymbol{T}^*, \boldsymbol{T}\big ]\!\!\big ]\big)$ is in trace class whenever $\boldsymbol{T}$ is in $BS_{m, \vartheta}(\Omega)$. 
\begin{thm}{\label{lem63}}
Let $\boldsymbol{T}=(T_1,\ldots, T_d)$ be a commuting tuple of operators on a Hilbert space $\mathcal{H}$. Assume that  $\boldsymbol{T}$ is in the class $BS_{m, \vartheta}(\Omega)$. Then the determinant operator $\text{dEt}\,\big(\big[\!\! \big [\boldsymbol{T}^*, \boldsymbol{T}\big ]\!\! \big ]\big)$ is in trace-class and \[\text{trace}\,\big (\text{dEt}\,\big(\big[\!\! \big [\boldsymbol{T}^*, \boldsymbol{T}\big ]\!\!\big]\big)\big )\leq m\, \vartheta \,d!\prod_{i=1}^{d}\|T_i\|^2.\]
\end{thm}
\begin{proof}
   By hypothesis, $$\big \|P_N\big (\sum_{\eta \in \mathfrak S_{d}}\text{Sgn} (\eta)T_{\eta(1)}^*T_{\tau(1)}T_{\eta(2)}^*\ldots T_{\tau(d-1)} T_{\eta(d)}^*\big )P_N^{\perp}T_{\tau(d)}P_N\big \|\leq \frac{\vartheta}{\binom{N+d-1}{d-1}}\prod_{i=1}^d \big \|T_i \big \|^2 .$$ Thus combining this inequality with the one from Lemma \ref{L1.6}, we have 
\begin{align*}
    \big |\text{trace}\,\big (P_N\, \text{dEt}\,\big(\big[\!\!\big[\boldsymbol{T}^*, \boldsymbol{T}\big]\!\!\big]\,\big) P_N\big )\big| \leq m\, \vartheta\, d!\prod_{i=1}^{d}\big \|T_i \big \|^2.
\end{align*}
Since $\text{dEt}\,\big(\big[\!\! \big [\boldsymbol{T}^*, \boldsymbol{T}\big ]\!\!\big ]\big)$ is non-negative definite by assumption and the projections $P_N$  increase to $I$ in the strong operator topology, we obtain the inequality
\begin{align*}
    \text{trace}\,\big(\text{dEt}\,\big(\big[\!\!\big[\boldsymbol{T}^*, \boldsymbol{T}\big ]\!\!\big ]\big) \big) \leq  m\, \vartheta\, d!\prod_{i=1}^{d}\big \|T_i\big \|^2
\end{align*}
completing the proof. 
\end{proof}
\subsection{The tensor product model}
For $i=1,2,$ let  $\boldsymbol T^{(i)}=(T^{(i)}_1,\ldots, T^{(i)}_{d_i})$ be a $d_i$- tuple of commuting bounded operators.  Set 
\begin{align*}(\boldsymbol{T}^{(1)}\,\#\, \boldsymbol{T}^{(2)}) &:=
(\boldsymbol T^{(1)}\otimes \boldsymbol I, \boldsymbol I \otimes \boldsymbol T^{(2)})\\ 
&= (T^{(1)}_1 \otimes I, \ldots , T^{(1)}_{d_1} \otimes I, I\otimes T^{(2)}_1, \ldots I\otimes ,T^{(2)}_{d_2}) \end{align*}
This definition clearly extends to $d_i$- tuples of commuting operators, $i=1, \ldots , n$.

\begin{lem} The spectrum $\sigma\big ( \boldsymbol{T}^{(1)}\,\#\, \boldsymbol{T}^{(2)} \big )$ of the operator $\boldsymbol{T}^{(1)}\,\#\, \boldsymbol{T}^{(2)}$ is $\sigma \big (\boldsymbol T^{(1)}\big ) \times \sigma \big ( \boldsymbol T^{(2)}\big )$. Moreover, 
if the $d_i$- tuples $\boldsymbol T^{(i)}$, $i=1,2$, are $m_i$- cyclic, then the operator $\boldsymbol{T}^{(1)}\,\#\, \boldsymbol{T}^{(2)}$ is $m$- cyclic, where $m\leq m_1 m_2$. 
\end{lem}
\begin{proof}
The joint spectrum of $\boldsymbol{T}^{(1)}\,\#\, \boldsymbol{T}^{(2)}$ is explicitly given in \cite[Theorem 2.2]{Vasilescu}. 
If $\boldsymbol \xi_{\boldsymbol T^{(i)}}[m_i]$, $i=1,2,$ is the cyclic set for the $d_i$- tuple  $\boldsymbol T^{(i)}$, then the cyclic set of the operator $\boldsymbol{T}^{(1)}\,\#\, \boldsymbol{T}^{(2)}$ is clearly contained in the set of vectors $$\big\{x\otimes y\mid x\in \boldsymbol \xi_{\boldsymbol T^{(1)}}[m_1] ~\text{and}~ y \in \boldsymbol \xi_{\boldsymbol T^{(2)}}[m_2]\big \}.$$ Thus the claim that $m \leq m_1 m_2$ is verified. 
\end{proof}

We now obtain a trace inequality for the operator $\text{dEt}\big (\big[\!\! \big [ (\boldsymbol{T}^{(1)}\,\#\, \boldsymbol{T}^{(2)})^*, (\boldsymbol{T}^{(1)}\,\#\, \boldsymbol{T}^{(2)}) \big ]\!\!\big ]\big )$.
A similar inequality can be proved for $\boldsymbol T^{(1)}\# \cdots \# \boldsymbol T^{(n)}.$
\begin{thm} \label{tensorprod}
Assume  that $\boldsymbol{T^{(i)}}$ is in the class $BS_{m_i,1}(\Omega_i)$, $i=1,2$. Then the determinant operator  $$\text{dEt}\,\big(\big[\!\! \big [ (\boldsymbol{T}^{(1)}\,\#\, \boldsymbol{T}^{(2)})^*, (\boldsymbol{T}^{(1)}\,\#\, \boldsymbol{T}^{(2)}) \big ]\!\!\big ]\big)$$ is non-negative definite and \[\text{trace}\,\big(\text{dEt}\,\big(\big[\!\! \big [ (\boldsymbol{T}^{(1)}\,\#\, \boldsymbol{T}^{(2)})^*, (\boldsymbol{T}^{(1)}\,\#\, \boldsymbol{T}^{(2)}) \big ]\!\!\big ]\big) \big)\leq 2 d_{1}!d_{2}!m_1 m_2 \prod_{i=1}^{d_1}\big \|T^{(1)}_i\big \|^2 \prod_{i=1}^{d_2}\big \|T^{(2)}_i\big \|^2. \]
\end{thm}
\begin{proof}
It is easy to see that 
\begin{align*}
\big[\!\! \big [ (\boldsymbol{T}^{(1)}\,\#\, \boldsymbol{T}^{(2)})^*, (\boldsymbol{T^{(1)}}\,\#\, \boldsymbol{T}^{(2)}) \big ]\!\!\big ]
= \begin{pmatrix}
\big[\!\! \big [ (\boldsymbol{T}^{(1)})^*, \boldsymbol{T}^{(1)} \big ]\!\!\big ]\otimes \boldsymbol{I} & 0\\
                   \\
0 & \boldsymbol{I}\otimes \big[\!\! \big [  (\boldsymbol{T}^{(2)})^*,  \boldsymbol{T}^{(2)} \big ]\!\!\big ]
\end{pmatrix}.  
\end{align*}
Thus $$\text{dEt}\,\big(\big[\!\! \big [ (\boldsymbol{T^{(1)}}\,\#\, \boldsymbol{T^{(2)}})^*, (\boldsymbol{T^{(1)}}\,\#\, \boldsymbol{T^{(2)}}) \big ]\!\!\big ]\big)=2\,\text{dEt}\,\big(\big[\!\! \big [ \boldsymbol{T^{(1)}}^*, \boldsymbol{T^{(1)}} \big ]\!\!\big ]\big) \, \otimes \text{dEt}\,\big(\big[\!\! \big [ \boldsymbol{T^{(2)}}^*, \boldsymbol{T^{(2)}} \big ]\!\!\big ]\big).$$
Since for $i=1,2,$ $\boldsymbol{T^{(i)}}$ is in the class $BS_{m_i,1}(\Omega_i)$, $\text{dEt}\,\big(\big[\!\! \big [ (\boldsymbol{T}^{(i)})^*, \boldsymbol{T}^{(i)} \big ]\!\!\big ]\big)$ is non-negative definite and 
\[\text{trace}\,\big(\text{dEt}\,\big(\big[\!\! \big [ (\boldsymbol{T}^{(i)})^*, \boldsymbol{T}^{(i)} \big ]\!\!\big ]\big) \big)\leq  d_{i}!m_i  \prod_{j=1}^{d_i}\big \|T^{(i)}_j\big \|^2 . \] 
Hence,  $\text{dEt}\,\big(\big[\!\! \big [ (\boldsymbol{T}^{(1)}\,\#\, \boldsymbol{T}^{(2)})^*, (\boldsymbol{T}^{(1)}\,\#\, \boldsymbol{T}^{(2)}) \big ]\!\!\big ]\big)$ is non-negative definite and 
\begin{align*}
 \text{trace}\, \big( \text{dEt}\,\big(\big[\!\! \big [ (\boldsymbol{T}^{(1)}\,\#\, \boldsymbol{T}^{(2)})^*&, (\boldsymbol{T}^{(1)}\,\#\, \boldsymbol{T}^{(2)}) \big ]\!\!\big ]\big)\big)\\&= 2\,\text{trace}\,\big(\text{dEt}\,\big(\big[\! \big [ (\boldsymbol{T}^{(1)})^*, \boldsymbol{T}^{(1)} \big ]\!\big ]\big) \big)\text{trace}\,\big(\text{dEt}\,\big(\big[\!\! \big [ (\boldsymbol{T}^{(2)})^*, \boldsymbol{T}^{(2)} \big ]\!\!\big ]\big) \big)\\
 &\leq 2 d_{1}!m_1  \prod_{i=1}^{d_1}\big \|T^{(1)}_i\big \|^2 \, \cdot d_{2}! m_2 \prod_{i=1}^{d_2}\big \|T^{(2)}_i\big \|^2 \\
 &=2 d_{1}!d_{2}!m_1 m_2 \prod_{i=1}^{d_1}\big \|T^{(1)}_i\big \|^2 \prod_{i=1}^{d_2}\big \|T^{(2)}_i\big \|^2.\qedhere
\end{align*} \end{proof}
\begin{rem}\begin{enumerate}
\item[]
\item Let $\boldsymbol T^{(i)}$, $i=1,\ldots ,n$ be a set of $n$ commuting $d_i$- tuple of operators. A similar proof, as given above, shows that if $\boldsymbol T^{(i)}\in \text{BS}_{m_i,1}(\Omega_i)$, then 
$$\text{dEt}\,\big(\big[\!\! \big [ (\boldsymbol{T}^{(1)}\,\#\cdots \# \boldsymbol{T}^{(n)})^*, (\boldsymbol{T}^{(1)}\,\#\cdots\# \boldsymbol{T}^{(n)}) \big ]\!\!\big ]\big)$$ is non-negative definite and 
\[\text{trace}\,\big(\text{dEt}\,\big(\big[\!\! \big [ (\boldsymbol{T}^{(1)}\,\# \cdots\#\, \boldsymbol{T}^{(n)})^*, (\boldsymbol{T^{(1)}}\,\# \cdots \#\, \boldsymbol{T}^{(n)}) \big ]\!\!\big ]\big) \big)\leq n!\,{d_{1}!\cdots d_{n}!\,m_1\cdots  m_n} \prod_{i=1}^{n}\|\boldsymbol T^{(i)}\|^2, \]
where $\|\boldsymbol T^{(i)}\|^2=\prod_{j=1}^{d_i}\|T^{(i)}_j\|^2$.
\item If $d_i=1$, $i=1,\ldots , n$, then  $(\boldsymbol T^{(1)} \# \cdots \# \boldsymbol T^{(n)})$ is of the form $(T_1 \otimes I \cdots \otimes I, \ldots I \otimes \cdots \otimes T_n)$. Now, we can apply the Berger-Shaw inequality to each of the operators $T_i$, $1\leq i \leq n$, to conclude 
$$\text{trace}\,\big(\text{dEt}\,\big(\big[\!\! \big [ (\boldsymbol{T}^{(1)}\,\# \cdots\#\, \boldsymbol{T}^{(n)})^*, (\boldsymbol{T^{(1)}}\,\# \cdots \#\, \boldsymbol{T}^{(n)}) \big ]\!\!\big ]\big) \big)\leq n!\,m_1\cdots  m_n \frac{\nu(\Omega_1\times \cdots \times \Omega_n)}{\pi^n}. $$
\end{enumerate}
\end{rem}
Let $\boldsymbol{M}=(M_1, \ldots M_d)$ be the $d$- tuple of multiplication by the coordinate functions  on the Hardy space $H^2(\mathbb{D}^d)$. Clearly, 
$\boldsymbol M = M\#\cdots \# M$ where $M$ is the multiplication operator on $H^2(\mathbb{D})$.
\begin{cor}\label{qmn}
For the $d$- tuple $\boldsymbol{M}= M\#\cdots \# M$ on the Hardy space $H^2(\mathbb{D}^d)$, we have that the operator $\text{dEt}\,\big(\big[\!\! \big [ \boldsymbol{M}^*, \boldsymbol{M} \big ]\!\!\big]\big)$ is non-negative definite  and $$\text{trace}\,\big(\text{dEt}\,\big(\big[\!\! \big [ \boldsymbol{M}^*, \boldsymbol{M} \big ]\!\!\big ]\big )\big) = \text{trace}\,\big(\text{dEt}\,\big(\big[\!\! \big [ (M\#\cdots \# M)^*, ( M\#\cdots \# M) \big ]\!\!\big ]\big )\leq d!.$$ 
\end{cor}
\begin{rem}
A direct computation shows that The inequality of Corollary \ref{qmn} is actually an equality. Consequently, it follows that the inequality obtained in Theorem \ref{tensorprod} is sharp. Although, the $d$- tuple $\boldsymbol{M}= M\#\cdots \# M$ is not in $BS_{m,\vartheta}(\Omega)$ for any choice of $m,\theta$, never the less, the trace estimate of Theorem \ref{lem63} agrees with the one obtained in Corollary \ref{qmn}.
\end{rem}

\section{Examples of operators in the class $BS_{m, \vartheta}(\Omega)$}

We consider two sets of examples, the first is based on the Euclidean ball $\mathbb B_d \subset \mathbb C^d$ while the second set of examples comes from considering the ball $\mathbb B_{2,1}:=\{(z_1,z_2): |z_1|^2 +  |z_2| < 1\} \subseteq \mathbb C^2$.  

\subsection{The case of the Euclidean ball $BS_{1, 1}(\mathbb B_d)$}
In the following examples, we have taken $\Omega = \mathbb B_d$. 
Let $H^{(\lambda)}(\mathbb B_d)$ be the weighted Bergman spaces of the unit Euclidean ball $\mathbb B_d$ discussed in Section \ref{weighted}.

Let $\mathcal U(d)$ be the group of unitary linear transformations on $\mathbb C^d$, let $\boldsymbol T$ be a commuting $d$-tuple of bounded linear operators on $\mathcal H$ and finally, let $\mathcal U(\mathcal H)$ be the group of unitary linear  transformations on $\mathcal H$. Clearly, the group $\mathcal U(d)$ acts on any commuting $d$-tuple of operators $\boldsymbol T$, namely, 
\begin{equation}\label{Sphericaldefn}
U\cdot \boldsymbol T: = \Big (\sum_{j=1}^d U_{1j} T_j, \ldots , \sum_{j=1}^d U_{d j} T_j \Big ), \: U= \big (\!\!\big ( U_{ij} \big )\!\!\big ) \in \mathcal U(d).
\end{equation} 
The $d$-tuple $\boldsymbol T$ is  said to be \textit{spherical} if there is a map $\Gamma:\mathcal U(d) \to \mathcal U(\mathcal H)$ such that 
\begin{align}
\Gamma_U \boldsymbol  T \Gamma_U^* := (\Gamma_U T_1 \Gamma_U^*, \ldots , \Gamma_U T_d \Gamma_U^*)
& = U\cdot \boldsymbol T\:\: \mbox{\rm for all } U\in \mathcal U(d).
\end{align} 

The set of vectors    
$$e_{\boldsymbol \alpha}(\boldsymbol z) = \sqrt{\tfrac{(d)_{|\boldsymbol \alpha|}}{\boldsymbol \alpha!}} z_1^{\alpha_1} \cdots z_d^{\alpha_d},\,\, \boldsymbol \alpha \in \mathbb N_0^d, \boldsymbol \alpha!=\alpha_1!\cdots \alpha_d!$$ is an orthonormal basis of $H^2(\mathbb{B}_d)$. 
The $d$- tuple $\boldsymbol S$ of multiplication operators  by the coordinate functions and its adjoint $\boldsymbol S^*$ on the Hardy space $H^2(\mathbb B_d)$ are commuting tuples of weighted shift operators: 
$$S_i e_{\boldsymbol \alpha} = \sqrt{\tfrac{\alpha_i +1 }{|\boldsymbol \alpha|+d}}\, e_{\boldsymbol \alpha +\epsilon_i},\,\, S^*_i e_{\boldsymbol \alpha} = \begin{cases}\sqrt{\tfrac{\alpha_i }{|\boldsymbol  \alpha|+d -1}}\, e_{\boldsymbol \alpha - \epsilon_i} & \text{\rm if}\,\, \alpha_i > 0, \\0 & \text{\rm otherwise}. \end{cases}$$
Note that the operator $\boldsymbol S$ is the same as the commuting tuple $\boldsymbol M^{(d)}$. However, it is convenient to use a different notation for this particular $d$-tuple as will be apparent soon. The basic properties of commuting tuples of weighted shifts, also called joint weighted shifts, are in \cite{JewelLubin}.


\begin{thm}\label{ThmI}
For the $d$ - tuple $\boldsymbol{S}$ of multiplication by the coordinate functions on the Hardy space $H^2(\mathbb B_d)$, the operator  $\text{dEt}\,\big (\big[\!\! \big [\boldsymbol{S}^*, \boldsymbol{S}\big ]\!\!\big ]\big)$ is non-negative definite and $\text{trace}\,\big (\text{dEt}\,\big (\big[\!\! \big [\boldsymbol{S}^*, \boldsymbol{S}\big ]\!\!\big ]\big)\big)=1.$
\end{thm}
\begin{proof}
The constant function $\boldsymbol 1$ is a cyclic vector for the $d$- tuple $\boldsymbol{S}$.  Also,  $\sigma(\boldsymbol{S})= \overline{\mathbb{B}}_d$, see \cite{CY}. To show that  $\text{dEt}\,\big (\big[\!\! \big [\boldsymbol{S}^*, \boldsymbol{S}\big ]\!\!\big ]\big)$ is non-negative definite,  we claim that 
\begin{align}\label{5.5}
\big( \sum_{\eta \in \mathfrak S_{d}} \text{sgn}&(\eta)S_{\eta(1)}^*S_{\tau(1)}\ldots S_{\tau(d-1)} S_{\eta(d)}^*\big)e_{\boldsymbol \alpha} \nonumber \\&= \text{sgn}(\tau)\sqrt{\frac{\alpha_{\tau(d)}}{|\boldsymbol \alpha|+d-1}}{(|\boldsymbol \alpha|+d-1)}^{-(d-1)}e_{\boldsymbol \alpha-\epsilon_{\tau(d)}}, \end{align}
for each fixed $\tau$ in $\mathfrak S_d$.
Here we have assumed $\alpha_{\tau(d)} > 0$ without loss of generality. 


For any fixed but arbitrary $\tau \in \mathfrak S_{k}$ and an arbitrary $k$- tuple $(i_1,\ldots ,i_{k})$, $k  \leq d$, with 
$1\leq i_1 < i_2 < \cdots < i_{k}  \leq d$, let ``$P_\tau(i_1, \ldots, i_k)$'' be the induction hypothesis, namely, the statement  
\begin{align}
\big( \sum_{\eta \in \mathfrak S_{d-1}}\text{sgn}&(\eta)S_{i_{\eta(1)}}^*S_{i_{\tau(1)}}\ldots S_{i_{\tau(k-1)}} S_{i_{\eta(k)}}^*\big)e_{\boldsymbol \alpha} \nonumber \\&= \text{sgn}(\tau)\sqrt{\frac{\alpha_{i_{\tau(k)}}}{|\boldsymbol \alpha|+d-1}}{(|\boldsymbol \alpha|+d-1)}^{-(k-1)}e_{\boldsymbol \alpha-\epsilon_{i_{\tau(k)}}}.\end{align}
We see that the equality of Equation \eqref{5.5} is the same as the equality asserted in  the statement $P_\tau(i_1, \ldots, i_d)$ with $k=d$ and $\tau = \text{id}$, the identity permutation on $d$ elements. 

Thus, a proof of the equality \eqref{5.5} follows from showing that  
the validity of all the equalities in $P_\tau(i_1, \ldots ,i_k)$,  $\tau\in \mathfrak{S}_k$,  $k < d$, implies the validity of all the equalities in $P_\tau(i_1, \ldots ,i_{k+1})$  $\tau\in \mathfrak{S}_{k+1}$.

To establish this, first note that if $k=1,$ then $P_\tau(i_1)$, $\tau \in S_1$ is the equalty ($\alpha_i > 0$):
$S^*_{i_1} e_{\boldsymbol \alpha} = \sqrt{\tfrac{\alpha_{i_1} }{|\boldsymbol \alpha|+d -1}}\, e_{\boldsymbol \alpha - \epsilon_{i_1}}$, which is clearly valid.  Now if $k=2$, we note that 
$$\big( \sum_{\eta \in \mathfrak S_{2}}\text{sgn}(\eta)S_{i_{\eta(1)}}^*S_{i_{\tau(1)}}S_{i_{\eta(2)}}^*\big)e_{\boldsymbol \alpha} \nonumber = \text{sgn}(\tau)\sqrt{\frac{\alpha_{i_{\tau(2)}}}{|\boldsymbol \alpha|+d-1}}{(|\boldsymbol \alpha|+d-1)}^{-1}e_{\boldsymbol \alpha-\epsilon_{i_{\tau(2)}}}$$
for any pair $i_1,i_2$ with $1\leq i_1 < i_2  \leq d$ and any fixed but arbitrary $\tau \in {\mathfrak S}_2$. This establishes the validity of $P_\tau(i_1,i_2)$, $\tau \in {\mathfrak S}_2$.

More generally, assume that the equalities in $P_\tau(i_1, \ldots , i_{k-1})$ are valid for every $(i_1,\ldots ,i_{k-1})$ with $1\leq i_1 < i_2  \cdots < i_{k-1}  \leq d$, and any fixed but arbitrary  $\tau \in \mathfrak S_{k-1}$. 
To show that the equality in $P_\tau(i_1, \ldots , i_{k})$ is valid 
for each fixed but arbitrary $(i_1, \ldots , i_k)$ and $\tau \in {\mathfrak S}_k$, it is enough to verify it with $\tau =id$. For this, we split the left hand side of $P_{id}(i_1, \ldots ,i_k)$ into several sums fixing $\eta(k)= j$, $j=1,\ldots ,k,$ in each one of these sums, that is, 
\begin{align}\label{twoterms}
\bigg( \sum_{\eta \in \mathfrak  S_{k}}&\text{sgn}(\eta)S_{i_{\eta(1)}}^*S_{i_1}S_{i_{\eta(2)}}^*\ldots S_{i_{\eta(k-1)}}^*S_{i_{k-1}} S_{i_{\eta(k)}}^*\bigg)e_{\boldsymbol \alpha}= \nonumber \\& \big(\sum_{\eta \in \mathfrak S_{k}, \eta(k) =k}\text{sgn}(\eta)S_{i_{\eta(1)}}^*S_{i_1}S_{i_{\eta(2)}}^*\ldots S_{i_{\eta(k-1)}}^*S_{i_{k-1}} S_{i_k}^*\big)e_{\boldsymbol \alpha} \nonumber \\&+\big(\sum_{\eta \in \mathfrak S_{k},\eta(k)=k-1}\text{sgn}(\eta)S_{i_{\eta(1)}}^*S_{i_1}S_{i_{\eta(2)}}^*\ldots S_{i_{\eta(k-1)}}^*S_{i_{k-1}} S_{i_{k-1}}^*\big)e_{\boldsymbol \alpha} \nonumber\\
&+\big(\sum_{\eta \in \mathfrak S_{k},\eta(k)=k-2}\text{sgn}(\eta)S_{i_{\eta(1)}}^*S_{i_1}S_{i_{\eta(2)}}^*\ldots S_{i_{\eta(k-1)}}^*S_{i_{k-1}} S_{i_{k-2}}^*\big)e_{\boldsymbol \alpha} \nonumber\\&+ \ldots +\big(\sum_{\eta \in \mathfrak S_{k}, \eta(k)=1}\text{sgn}(\eta)S_{i_{\eta(1)}}^*S_{i_1}S_{i_{\eta(2)}}^*\ldots S_{i_{\eta(k-1)}}^*S_{i_{k-1}} S_{i_1}^*\big)e_{\boldsymbol \alpha}.
\end{align}
Pick a fixed but arbitrary sum in $P_{id}(i_1, \ldots ,i_k)$ 
with $\eta(k) = j$, $j=1, \ldots, k-2$. We claim that these sums vanish.   
Each one of these sums is of the form
\begin{equation}\label{oneterm}
\sum_{\eta\in \mathfrak  S_k, \eta(k)=j} \text{sgn}(\eta)S_{i_{\eta(1)}}^*S_{i_1}S_{i_{\eta(2)}}^*\ldots S_{i_{j-1}} S^*_{i_{\eta(j)}}S_{i_j}S^*_{i_{\eta(j+1)}} S_{i_{j+1}}S^*_{i_{\eta(j+2)}}
\ldots S_{i_{\eta(k-1)}}^*S_{i_{k-1}} S_{i_j}^*.
\end{equation}
For a fixed $\eta\in {\mathfrak S}_k$, let $\eta_\sigma\in {\mathfrak S}_k$ be the permutation:
$$\eta_\sigma(i) = \begin{cases}\eta(i) & i\not \in \{ j, j+1\},\\
\eta(j+1) & \text{\rm if}\,\, i = j, \\
\eta(j) & \text{\rm if}\,\, i= j+1. \end{cases}$$
The sign of $\eta_\sigma$ is opposite of the sign of $\eta$ and these occur in pairs. Also, $S_i^*S_l S^*_p = S^*_p S_l S_i^*$ for any choice of ($(i,l,p)$ with $i\not = l \not = p$. Clearly, $\eta(j)\not = j \not = \eta(j+1)$ by choice. 
Putting these observations together, we conclude that the sum \eqref{oneterm} vanishes.

Now, we examine the two nonzero sums that remain on the right hand side of \eqref{twoterms}. The first of these two sums is
\begin{align*}
    \big(\sum_{\eta \in \mathfrak S_{k}}&\text{sgn}(\eta)S_{i_{\eta(1)}}^*S_{i_1}S_{i_{\eta(2)}}^*\ldots S_{i_{\eta(k-1)}}^*S_{i_{k-1}} S_{i_k}^*\big)e_{\boldsymbol \alpha}\\=&  \big(\sum_{\sigma \in \mathfrak S_{k-1}}\text{sgn}(\sigma)S_{i_{\sigma(1)}}^*S_{i_1}S_{i_{\sigma(2)}}^*\ldots S_{i_{\sigma(k-1)}}^*\big)S_{i_{k-1}} S_{i_k}^*e_{\boldsymbol \alpha}.
\end{align*}
Applying the equality in $P_{id}(i_1, \ldots, i_{k-1})$ to the vector $S_{i_{k-1}} S_{i_k}^*e_{\boldsymbol \alpha}$, we have 
\begin{equation}\label{First}
\big(\sum_{\eta \in \mathfrak S_{k}}\text{sgn}(\eta)S_{i_{\eta(1)}}^*S_{i_1}S_{i_{\eta(2)}}^*\ldots S_{i_{\eta(k-1)}}^*S_{i_{k-1}} S_{i_k}^*\big)e_{\boldsymbol \alpha} = \sqrt{\frac{\alpha_{i_k}}{|\boldsymbol \alpha|+d-1}}\frac{\alpha_{i_{k-1}}+1}{(|\boldsymbol \alpha|+d-1)^{(k-1)}}e_{\boldsymbol \alpha-\epsilon_{i_k}}.
\end{equation}
The second sum 
\begin{align*}
\big(\sum_{\eta \in \mathfrak S_{k}}&\text{sgn}(\eta)S_{i_{\eta(1)}}^*S_{i_1}S_{i_{\eta(2)}}^*\ldots S_{i_{\eta(k-1)}}^*S_{i_{k-1}} S_{i_{k-1}}^*\big)e_{\boldsymbol \alpha}\\=&
-\big(\sum_{\sigma \in \mathfrak S_{k-1}}\text{sgn}(\sigma)S_{i_{\sigma(1)}}^*S_{i_1}S_{i_{\sigma(2)}}^*\ldots S_{i_{\sigma(k)}}^*\big)S_{i_{k-1}} S_{i_{k-1}}^*e_{\boldsymbol \alpha}.
\end{align*}
The verification of this equality is an immediate consequence of the following observations. 
\begin{enumerate} 
\item Each permutation $\eta \in \hat{\mathfrak S}_k:=\{\eta\in \mathfrak S_k: \eta(k) =k-1 \}$ is obtained by extending
some bijection $\hat{\sigma}: \{1,\ldots , k-1\} \to \{1, \ldots , k-2, k\}$ to a bijection of the set $\{1, \ldots ,k\}$.
\item The bijections $\hat{\sigma}$ are in one to one correspondence with the group of  permutations on a set of $k-1$ elements. For convenience, we take this set to be $\mathbb E_{k-1}:=\{1, \ldots , k-2, k\}$. 
\item Since $\eta\in \hat{\mathfrak S}_k$ is an extension of a permutation $\hat{\sigma}$ on $\mathbb E_{k-1}$ obtained by setting $\eta=\hat{\sigma}$ on $\mathbb E_{k-1}$ ( $\eta(k)=k-1$), it follows that $\text{sgn}(\eta) = -\text{sgn}(\hat{\sigma})$.
\item The group of permutations on the set $\mathbb E_{k-1}$ is isomorphic to $\hat{\mathfrak S}_{k-1}$ via the map $\hat{\sigma }\to \eta$.
\end{enumerate}
Now, as before, applying the equality in $P_{id}(i_1, \ldots, i_{k-2}, i_k)$ to the vector $S_{i_{k-1}} S_{i_{k-1}}^*e_{\alpha}$, we have  
\begin{align}\label{Second}
\big(\sum_{\eta \in \mathfrak S_{k}} \text{sgn}(\eta)S_{i_{\eta(1)}}^*S_{i_1}S_{i_{\eta(2)}}^*\ldots S_{i_{\eta(k-1)}}^*S_{i_{k-1}} S_{i_{k-1}}^*\big)e_{\boldsymbol \alpha} =&
- \sqrt{\frac{\alpha_{i_k}}{|\boldsymbol \alpha|+d-1}}\frac{\alpha_{i_{k-1}}}{(|\boldsymbol \alpha|+d-1)^{(k-1)}}e_{\boldsymbol \alpha-\epsilon_{i_k}}.
\end{align}
Combining \eqref{First} and \eqref{Second}, we obtain 
\begin{align*}
\big( \sum_{\eta \in \mathfrak S_{k}}\text{sgn}(\eta)S_{i_{\eta(1)}}^*S_{i_1}S_{i_{\eta(2)}}^*\ldots S_{i_{\eta(k-1)}}^*S_{i_{k-1}} S_{i_{\eta(k)}}^*\big)e_{\boldsymbol \alpha}=& \sqrt{\frac{\alpha_{i_k}}{|\boldsymbol \alpha|+d-1}}(|\boldsymbol \alpha|+d-1)^{-(k-1)}e_{\boldsymbol \alpha-\epsilon_{i_k}},
\end{align*}
which verifies $P_\tau(i_1, \ldots i_k)$ for $\tau=id$.  The verification of \eqref{5.5} for an arbitrary choice of $\tau \in {\mathfrak S}_k$ is similar.
Recall the expression of determinant operator in \eqref{det}
\begin{align*}
     \text{dEt}\big(\big[\!\! \big [\boldsymbol{S}^*, \boldsymbol{S}\big ]\!\!\big ]\big)=&\sum_{\tau,\eta\in \mathfrak  S_d}\text{sgn}(\tau)Sgn(\eta)S^*_{\eta(1)}S_{\tau(1)}S^*_{\eta(2)}\ldots S^*_{\eta(d)}S_{\tau(d)}+\\& -\sum_{\tau,\eta\in \mathfrak S_d}\text{sgn}(\tau)\text{sgn}(\eta)S_{\tau(d)}S^*_{\eta(1)}S_{\tau(1)}\ldots S_{\tau(d-1)}S_{\eta(d)}^*.
\end{align*}
Now using the equality in \eqref{5.5} we get,
\begin{align*}
    \sum_{\tau,\eta\in \mathfrak S_d}\text{sgn}(\tau)\text{sgn}(\eta)S^*_{\eta(1)}S_{\tau(1)}S^*_{\eta(2)}\ldots S^*_{\eta(d)}S_{\tau(d)}e_{\boldsymbol \alpha}=\frac{(d-1)!}{(|\boldsymbol \alpha|+d)^{(d-1)}}e_{\boldsymbol \alpha},
\end{align*}
and
\begin{align*}
    \sum_{\tau,\eta\in \mathfrak S_d}\text{sgn}(\tau)\text{sgn}(\eta)S_{\tau(d)}S^*_{\eta(1)}S_{\tau(1)}\ldots S_{\tau(d-1)}S_{\eta(d)}^*e_{\boldsymbol \alpha}=\frac{(d-1)!|\boldsymbol \alpha|}{(|\boldsymbol \alpha|+d-1)^{d}} e_{\boldsymbol \alpha}.
\end{align*}
Thus combining the above equalities we get
\begin{align*}
    \text{dEt}\,\big (\big[\!\! \big [\boldsymbol{S}^*, \boldsymbol{S}\big ]\!\!\big ]\big)e_{\boldsymbol \alpha}=\big( \frac{(d-1)!}{(|\boldsymbol \alpha|+d)^{(d-1)}}-\frac{(d-1)!|\boldsymbol \alpha|}{(|\boldsymbol \alpha|+d-1)^{d}} \big )e_{\boldsymbol \alpha}.
\end{align*}
It is easy to see $ \text{dEt}\,\big (\big[\!\! \big [\boldsymbol{S}^*, \boldsymbol{S}\big ]\!\!\big ]\big)$ is non-negative definite since $$\frac{1}{(|\boldsymbol \alpha|+d)^{(d-1)}}-\frac{|\boldsymbol \alpha|}{(|\boldsymbol \alpha|+d-1)^{d}}\geq 0.$$
Now, \begin{align*}
    \text{trace}\,\big (\text{dEt}\,\big (\big[\!\! \big [\boldsymbol{S}^*, \boldsymbol{S}\big ]\!\!\big ]\big)\big)=&\sum_{\alpha \in \mathbb{N}_0^d} \big \langle \text{dEt}\,\big (\big[\!\! \big [\boldsymbol{S}^*, \boldsymbol{S}\big ]\!\!\big ]\big)e_{\boldsymbol \alpha}, e_{\boldsymbol \alpha} \big \rangle\\
    =& \sum_{\boldsymbol \alpha \in \mathbb{N}_0^d}(d-1)!\big(\frac{1}{(|\boldsymbol \alpha|+d)^{(d-1)}}-\frac{|\boldsymbol \alpha|}{(|\boldsymbol \alpha|+d-1)^{d}} \big)\\
    =& \sum_{k=0}^{\infty}\sum_{\substack{\alpha_1,\ldots, \alpha_d \\ |\boldsymbol\alpha|=k}}(d-1)!\big(\frac{1}{(|\boldsymbol \alpha|+d)^{(d-1)}}-\frac{|\boldsymbol \alpha|}{(|\boldsymbol \alpha|+d-1)^{d}} \big)\\
    =& \sum_{k=0}^{\infty}\bigg(\frac{(k+d-1)(k+d-2)\ldots(k+1)}{(k+d)^{(d-1)}}-\frac{(k+d-2)\ldots(k+1) k}{(k+d-1)^{d-1}} \bigg)\\
    =& 1.
   \end{align*}
     This completes the proof.
\end{proof}
For $\boldsymbol a\in \mathbb C^d$ and $r >0$, let $\mathbb B[\boldsymbol a ,r]$
be the ball $\{\boldsymbol z\in \mathbb C^d: \|\boldsymbol z - \boldsymbol a\|_2 < r\}.$ We let $\mathbb B[r]$ denote the ball of radius $r$ centred at $0$. Finally, $\mathbb B_d$ is the unit ball in $\mathbb C^d$.

\begin{thm}\label{ThmII}
Let $\boldsymbol T$ be a $d$- tuple of spherical joint weighted shift operators and $T_{\delta}$ be the one variable weighted shift corresponding to $\boldsymbol{T}$. If $T_\delta$ is hyponormal, then $\boldsymbol T$ is in $BS_{1,1}(\mathbb{B}[r])$, where $\mathbb B[r] = \{\boldsymbol z\in \mathbb C^d : \|\boldsymbol z\|_2 < r\}$, $r >0$. \end{thm}

\begin{proof} As in the case of two variables, the $d$- tuple $\boldsymbol T$ is $1$- cyclic, $P_N T_j P_N^{\perp}=0,\,1\leq j\leq d,$ and its spectrum is of the form $\mathbb B[r]$ for some $r > 0$ depending on $\{\delta_{|\boldsymbol \alpha|}\}$. 

A calculation similar to the one in the proof of Theorem
\ref{ThmI} given above shows that 
\begin{align}\label{computation}
    \text{dEt}\,\big (\big[\!\! \big [\boldsymbol{T}^*, \boldsymbol{T}\big ]\!\!\big ]\big)x_{\boldsymbol \alpha}=\big( \frac{(d-1)!\delta_{|\boldsymbol \alpha|}^{2d}}{(|\boldsymbol \alpha|+d)^{(d-1)}}-\frac{(d-1)!|\boldsymbol \alpha|\delta_{|\boldsymbol \alpha|-1}^{2d}}{(|\boldsymbol \alpha|+d-1)^{d}} \big )x_{\boldsymbol \alpha}.
\end{align}
Since $\delta_{|\boldsymbol \alpha|}$ is an increasing sequence, it follows $$\big( \frac{(d-1)!\delta_{|\boldsymbol \alpha|}^{2d}}{(|\boldsymbol \alpha|+d)^{(d-1)}}-\frac{(d-1)!|\boldsymbol \alpha|\delta_{|\boldsymbol \alpha|-1}^{2d}}{(|\boldsymbol \alpha|+d-1)^{d}} \big )\geq (d-1)! \delta_{|\boldsymbol \alpha|-1}^{2d}\big( \frac{1}{(|\boldsymbol \alpha|+d)^{(d-1)}}-\frac{|\boldsymbol \alpha|}{(|\boldsymbol \alpha|+d-1)^{d}} \big ) \geq 0.$$ Hence $\text{dEt}\,\big (\big[\!\! \big [\boldsymbol{T}^*, \boldsymbol{T}\big ]\!\! \big ]\big)$ is non-negative definite.

To complete the proof, we need to verify the norm estimate (iii) of Definition \ref{class}. For this, taking $\tau$ to be the identity permutation in \eqref{5.5}, we obtain the equality  
\begin{align} \label{tau=id}
\big( \sum_{\eta \in \mathfrak S_{d}} \text{sgn}&(\eta)T_{\eta(1)}^*T_{1}T_{\eta(2)}^*\ldots T_{d-1} T_{\eta(d)}^*\big)x_{\boldsymbol \alpha} \nonumber \\&= \delta_{|\boldsymbol \alpha|}^{2d-1}\sqrt{\frac{\alpha_{d}}{|\boldsymbol \alpha|+d-1}}{(|\boldsymbol \alpha|+d-1)}^{-(d-1)}x_{\boldsymbol \alpha-\epsilon_{d}}.
\end{align}
Clearly, we have 
\begin{align}\label{final ineq}
    \|P_N \big( &\sum_{\eta \in \mathfrak S_{d}} \text{sgn}(\eta)T_{\eta(1)}^*T_{1}T_{\eta(2)}^*\ldots T_{d-1} T_{\eta(d)}^*\big) P_N^{\perp}\|\nonumber \\
    &\leq {\binom{N+d-1}{d-1}}^{-1}\prod_{\{i:\,i\neq d\}} \|T_i\|^2 \|T_{d}\|. \end{align}
    Consequently, \begin{align}\label{second final ineq}
    \|P_N \big( &\sum_{\eta \in \mathfrak S_{d}} \text{sgn}(\eta)T_{\eta(1)}^*T_{1}T_{\eta(2)}^*\ldots T_{d-1} T_{\eta(d)}^*\big) P_N^{\perp}T_{d} P_N\|\leq {\binom{N+d-1}{d-1}}^{-1}\prod_{i=1}^d \|T_i\|^2. \end{align}
It follows from Remark \ref{general tau}(b) that the inequality in Equation \eqref{second final ineq} remains unchanged when we replace the identity permutation by any other permutation from $\mathfrak S_d$.  Therefore, the $d$- tuple $\boldsymbol T$ is in the class $BS_{1, 1}(\mathbb B[r])$ and the  proof is complete.
\end{proof}

\begin{cor}\label{ballvolume}
Suppose $\delta_{n}\uparrow 1$ then 
$\text{trace}\,\big(\text{dEt}\,(\big[\!\! \big [\boldsymbol{T}^*, \boldsymbol{T}\big ]\!\!\big ])\big)=1.$
\end{cor}
\begin{proof}
The  string of equalities 
\begin{align*}
    \text{trace}\,\big (\text{dEt}\,\big (\big[\!\! \big [\boldsymbol{T}^*, \boldsymbol{T}\big ]\!\!\big ]\big)\big)=&\sum_{\alpha \in \mathbb{N}^d} \big \langle \text{dEt}\,\big (\big[\!\! \big [\boldsymbol{T}^*, \boldsymbol{T}\big ]\!\!\big ]\big)x_{\boldsymbol \alpha}, x_{\boldsymbol \alpha} \big \rangle\\
    =& \sum_{\boldsymbol \alpha \in \mathbb{N}_0^d}(d-1)!\big(\frac{\delta_{|\alpha|}^{2d}}{(|\boldsymbol \alpha|+d)^{(d-1)}}-\frac{\delta_{|\alpha|-1}^{2d}|\boldsymbol \alpha|}{(|\boldsymbol \alpha|+d-1)^{d}} \big)\\
    =& \sum_{k=0}^{\infty}\sum_{\substack{\alpha_1,\ldots, \alpha_d \\ |\boldsymbol\alpha|=k}}(d-1)!\big(\frac{\delta_{|\alpha|}^{2d}}{(|\boldsymbol \alpha|+d)^{(d-1)}}-\frac{\delta_{|\alpha|-1}^{2d}|\boldsymbol \alpha|}{(|\boldsymbol \alpha|+d-1)^{d}} \big)\\
    =& \sum_{k=0}^{\infty}\bigg(\delta_k^{2d}\frac{(k+d-1)(k+d-2)\ldots(k+1)}{(k+d)^{(d-1)}}-\delta_{k-1}^{2d}\frac{(k+d-2)\ldots(k+1) k}{(k+d-1)^{d-1}} \bigg)\\
    =&\lim_{k\rightarrow \infty}\delta_k^{2d}\frac{(k+d-1)(k+d-2)\ldots(k+1)}{(k+d)^{(d-1)}}=1.
\end{align*}
where second equality follows from Equation \eqref{computation}, verifies the claim.
\end{proof}

\subsection{The case of an ellipsoid $BS_{1, 2}(\mathbb{B}_{2, 1})$}
For $p,q\in \mathbb N$,  let
$\mathbb B_{p,q}=\big \{\boldsymbol z \in \mathbb{C}^2: |z_1|^p+|z_2|^q<1 \big\}$. These are examples of pseudo convex Reinhardt domains in $\mathbb C^2$. The usual Euclidean ball $\mathbb B_2$ is obtained by taking $p=q=2$, i.e., $\mathbb B_{2,2}= \mathbb B_2$. 

The pair $(z_1,z_2)\in \mathbb C^2$ is in $\mathbb B_{2,1}$ if and only if $r_1^2 + r_2 < 1,$ where $r_{k}:=|z_{k}|$, $k=1,2$. The volume measure $\nu$ restricted to $\mathbb B_{2,1}$ is of the form $d\nu(\boldsymbol z)=r_1 r_2 d r_1 d r_2 d\theta_1 d\theta_2$, $z_k= r_k\exp(i\theta_k)$, $k=1,2$ and set 
\begin{equation}
d\mu_{\lambda}(z):= (1-r_1^2-r_2)^{\lambda-4}r_1 r_2 d r_1 d r_2 d\theta_1 d\theta_2.  
\end{equation} 
The measure $d\mu_\lambda$ defines an inner product on the space $\mathbb C[\boldsymbol z]$ of polynomials in two variables by integration over $\mathbb B_{2,1}$:
$$\langle p , q \rangle_{\lambda} := \int_{\mathbb B_{2,1}} p \overbar{q} d \mu_\lambda.$$ Let $\mathbb A^{(\lambda)}(\,\mathbb B_{2,1})$
denote the Hilbert space obtained by taking the completion of the inner product space $\big (\mathbb C[\boldsymbol z], \langle \cdot, \cdot \rangle_\lambda\big )$. The Hilbert space $\mathbb A^{(\lambda)}(\,\mathbb B_{2,1})$ is non-zero if and only if $\lambda > 3$. This follows from the norm computation below. 
For any multi-index $\boldsymbol \alpha=(\alpha_1, \alpha_2) \in \mathbb N_0^2$, we have 
\begin{align*}
    \|\boldsymbol{z}^{\boldsymbol \alpha}\|_{\lambda}^2=& \int_{\Omega_{2,1}}|\boldsymbol {z}^{\boldsymbol \alpha}|^2 d{\mu}_{\lambda}\\
    =& (2\pi)^2 \int_{r_1=0}^{1}\int_{r_2=0}^{1-r_1^2}r_1^{2\alpha_1+1}r_2^{2\alpha_2 +1}(1-r_1^2-r_2)^{\lambda-4}d r_1 d r_2\\
    =& 2 (\pi)^2 B(2\alpha_2 +2, \lambda-3)B(\alpha_1 +1, 2\alpha_2 +\lambda -1)\\
    =& 2 (\pi)^2 \Gamma(\lambda-3)\frac{\Gamma(\alpha_1 +1)\Gamma(2\alpha_2 +2)}{\Gamma(2\alpha_2 +\alpha_1+\lambda)}.
    \end{align*}
Integrating first, with respect to the measure $d\theta_1 d\theta_2$, we see that $\{\boldsymbol z^{\boldsymbol \alpha} \mid \boldsymbol \alpha\in \mathbb N_0^2\}$ is an orthogonal set of vectors relative to the inner product $\langle \cdot, \cdot \rangle_\mu$ and hence 
the set of vectors $\big \{\phi_{\boldsymbol \alpha}:= \frac{\boldsymbol{z}^{\boldsymbol \alpha}}{\|\boldsymbol{z}^{\boldsymbol \alpha}\|_{\lambda}} : \boldsymbol \alpha \in \mathbb N_0^2 \big \}$ is a complete orthonormal set in the Hilbert space $\mathbb A^{(\lambda)}(\,\mathbb B_{2,1})$. Now, it is easy to see that 
the multiplication operators $M_{z_i}$, $i=1,2$ on the Hilbert space $\mathbb A^{(\lambda)}(\,\mathbb B_{2,1})$ are weighted shifts relative to this orthonormal basis, that is, 
$M_{z_i}(\phi_{\boldsymbol \alpha})=w_{\boldsymbol \alpha}^{(i)}\phi_{\boldsymbol \alpha+\epsilon_i}$, where the weights are given explicitly by the formulae:
$$w_{\boldsymbol \alpha}^{(1)}=\sqrt{\frac{\alpha_i+1}{\alpha_1+2\alpha_2+\lambda}}~\text{and}~w_{\boldsymbol \alpha}^{(2)}=\sqrt{\frac{(2\alpha_2+2)(2\alpha_2+3)}{(\alpha_1+2\alpha_2+\lambda)(\alpha_1+2\alpha_2+\lambda+1)}}.$$
Since  $\text{Sup}\{w_{\boldsymbol \alpha}^{(i)}=1\}$, it follows that $\|M_{z_i}\|=1$, $i=1,2.$ Many more details and the spectral picture of this pair of operators is given in \cite[Example 5.2]{Curto91spectral}
\begin{thm}
Let $\boldsymbol M=(M_{z_1}, M_{z_2})$ be the pair of  multiplication operators on $\mathbb A^{(\lambda)}(\,\mathbb B_{2,1})$ by the co-ordinate functions. If $\lambda\geq 4$, then $\boldsymbol M$ is in $BS_{1,2}(\,\mathbb{B}_{2,1})$. \end{thm}
\begin{proof}
Since $\boldsymbol M$ is a pair of joint weighted shifts, by definition, $P_N M_j P_N^{\perp}=0,\,i=1,2.$
The commuting pair $\boldsymbol M$ is $1$-cyclic 
and the Taylor joint spectrum $\sigma(\boldsymbol M) = \overbar{\mathbb B}_{2,1}$, see \cite{Curto91spectral}. The following computation verifies the estimate (iii) of the Definition \ref{class}. 
\begin{align*}
    &\big(M_{z_1}^*M_{z_1}M_{z_2}^*-M_{z_2}^*M_{z_1}M_{z_1}^*\big)\phi_{\boldsymbol \alpha}\\&=\Big(\, \frac{(2\alpha_1+2\alpha_2+\lambda-1)}{(\alpha_1+2\alpha_2+\lambda-2)(\alpha_1+2\alpha_2+\lambda-1)}\sqrt{\frac{(2\alpha_2)(2\alpha_2+1)}{(\alpha_1+2\alpha_2+\lambda-2)(\alpha_1+2\alpha_2+\lambda-1)}}\,\Big )\,\phi_{\boldsymbol \alpha-\epsilon_2}.
\end{align*}
Thus $\|P_N\big(M_{z_1}^*M_{z_1}M_{z_2}^*-M_{z_2}^*M_{z_1}M_{z_1}^*)P_N^{\perp}\|\leq \frac{2}{N+1}\|M_{z_1}\|^2\|M_{z_2}\|$ and therefore 
$$\|P_N\big(M_{z_1}^*M_{z_1}M_{z_2}^*-M_{z_2}^*M_{z_1}M_{z_1}^*)P_N^{\perp}M_{z_2}P_N\|\leq \frac{2}{N+1}\|M_{z_1}\|^2\|M_{z_2}\|^2.$$
\begin{multline*}
    \big (M_{z_2}^*M_{z_2}M_{z_1}^*-M_{z_1}^*M_{z_2}M_{z_2}^*\big)\phi_{\boldsymbol \alpha}\\=\Bigg(\frac{(2\alpha_2+1)(\alpha_1+\lambda-2)+2(\alpha_2+1)(\alpha_1+2\alpha_2+\lambda-2)}{(\alpha_1+2\alpha_2+\lambda-2)(\alpha_1+2\alpha_2+\lambda)}\\ \frac{2}{(\alpha_1+2\alpha_2+\lambda-1)}\sqrt{\frac{\alpha_1}{\alpha_1+2\alpha_2+\lambda-1}}\,\Bigg )\,\phi_{\boldsymbol \alpha-\epsilon_1}. 
\end{multline*}
Consequently, we have  $$\big \|P_N\big(M_{z_2}^*M_{z_2}M_{z_1}^*-M_{z_1}^*M_{z_2}M_{z_2}^*)P_N^{\perp}\big \|\leq \frac{2}{N+1}\|M_{z_2}\|^2\|M_{z_1}\|$$ and
$$\big \|P_N\big(M_{z_2}^*M_{z_2}M_{z_1}^*-M_{z_1}^*M_{z_2}M_{z_2}^*)P_N^{\perp}M_{z_1}P_N \big \|\leq \frac{2}{N+1}\|M_{z_2}\|^2\|M_{z_1}\|^2.$$

To complete the proof, we only need to verify that the operator 
$\text{dEt}\,\big(\big[\!\! \big [\boldsymbol{M}^*, \boldsymbol{M}\big ]\!\!\big ]\big)$ is non-negative definite. Evaluating on  the orthonormal basis $\{\phi_{\boldsymbol \alpha}\}$, we see that $\text{dEt}\,\big(\big[\!\! \big [\boldsymbol{M}^*, \boldsymbol{M}\big ]\!\!\big ]\big)\phi_{\boldsymbol \alpha} = \chi_{\boldsymbol \alpha} \phi_{\boldsymbol \alpha },$ where 
\begin{multline*}
\chi_{\boldsymbol \alpha} =  -\frac{2 \alpha_2 (2\alpha_2 +1) (2 \alpha_1+2 \alpha_2+\lambda-1)}{(\alpha_1+2 \alpha_2+\lambda-2)^2 (\alpha_1+2 \alpha_2+\lambda-1)^2}+\frac{(2\alpha_2 +2) (2 \alpha_2+3) (2 \alpha_1+2 \alpha_2+\lambda+1)}{(\alpha_1+2 \alpha_2+\lambda)^2 (\alpha_1+2 \alpha_2+\lambda+1)^2}\\+\frac{2 (\alpha_1+1) ((2 \alpha_2+1) (\alpha_1+\lambda-1)+2 (\alpha_2+1) (\alpha_1+2 \alpha_2+\lambda-1))}{(\alpha_1+2\alpha_2 +\lambda-1) (\alpha_1+2 \alpha_2+\lambda)^2 (\alpha_1+2 \alpha_2+\lambda+1)}\\-\frac{2 \alpha_1 ((2 \alpha_2+1) (\alpha_1+\lambda-2)+2 (\alpha_2+1) (\alpha_1+2 \alpha_2+\lambda-2))}{(\alpha_1+2 \alpha_2+\lambda-2) (\alpha_1+2 \alpha_2+\lambda-1)^2 (\alpha_1+2 \alpha_2+\lambda)}. \end{multline*}
Gathering these terms over a common denominator and simplifying we find that $\chi_{\boldsymbol \alpha}$ is a fraction with a positive denominator and the numerator is $4$ times the expression given below: 
\begin{multline*}
\!\!\!\! \alpha_1^4 (4 \alpha_2 \lambda-10 \alpha_2+3 \lambda-9)+2 \alpha_1 \left(8 \alpha_2^3 \left(12 \lambda^2-30 \lambda+1\right)+\alpha_2^2 \left(48 \lambda^3-96 \lambda^2-84 \lambda+98\right  )  \right ) \\ +\,2 \alpha_1^3 \left(8 \alpha_2^2 (2 \lambda-5)+2 \alpha_2 \left(4 \lambda^2-6 \lambda-7\right)+3 \left(2 \lambda^2-7 \lambda+6\right)  \right)+\alpha_1^2 \left(4 \alpha_2^2 \left(24 \lambda^2-54 \lambda-5 \right ) \right ) \\
+\,\alpha_1^2 \left ( 48 \alpha_2^3 (2 \lambda-5)+4 \alpha_2 \left(6 \lambda^3-3 \lambda^2-35 \lambda+32\right)  +18 \lambda^3-72 \lambda^2+93 \lambda-33\right) \\ 
+\, 2 \alpha_1 \left ( 32 \alpha_2^4 (2 \lambda-5)+\alpha_2 \left(8 \lambda^4+4 \lambda^3-98 \lambda^2+128 \lambda-39\right)+3 \left(2 \lambda^4-9 \lambda^3+13 \lambda^2-5 \lambda-2\right)\right) + \\ 
+\,   16 \alpha_2^4 \left(8 \lambda^2-21 \lambda+2\right)+16 \alpha_2^3 \left(6 \lambda^3-15 \lambda^2-4 \lambda+9\right)+4 \alpha_2^2 \left(8 \lambda^4-14 \lambda^3-37 \lambda^2+61 \lambda-14\right) \\
+\,  32 \alpha_2^5 (2 \lambda-5) +2 \alpha_2 \left(2 \lambda^5+3 \lambda^4-42 \lambda^3+64 \lambda^2-14 \lambda-13\right)+3 (\lambda+1) \left(\lambda^2-3 \lambda+2\right)^2.
\end{multline*}
It is then not hard to verify that the the constant term and the coefficients of $\alpha_1^i$, $1\leq i \leq 4$ together with that of $\alpha_2^j$, $1\leq j \leq 5$,   are all positive if $\lambda \geq 4$ completing the proof. 
\end{proof}

\begin{rem}
Although, unlike the case of the Euclidean ball $\mathbb B_2$, we are not able to explicitly compute the trace of the operator ${\rm dEt}\big(\big[\!\! \big [ \boldsymbol T^* , \boldsymbol T\big ]\!\!\big]\big)$ for the weighted Bergman spaces $\mathbb A^{(\lambda)}(\mathbb B_{2,1})$, extensive numerical computations show that it is approximately equal to $\tfrac{2}{3}$, which is the $\tfrac{2}{\pi^2}$ times the volume of the ellipsoid $\mathbb B_{2,1}$. 
\end{rem}
From Theorems 7.1 and 7.2 of \cite{HeltonHowe}, it follows that the trace of the generalized commutator $\text{\rm GC}(\boldsymbol T^*,\boldsymbol T)$ of a class of analytic Toeplitz operators is bounded above by $1$. In particular, the explicit formula given in Theorem 7.2 (a) of \cite{HeltonHowe} shows that equality is achieved for the tuple of multiplication by the coordinate functions. Of course, the same is true of the determinant operator ${\rm dEt}\big(\big[\!\! \big [ \boldsymbol T^* , \boldsymbol T\big ]\!\!\big]\big)$. In the example of weighted Bergman spaces over the Euclidean ball $\mathbb B_{d}$, from Corollary \ref{ballvolume}, it follows that 
$$\text{\rm trace}\,\big (\text{dEt}\,\big(\big[\!\! \big [\boldsymbol{T}^*, \boldsymbol{T}\big ]\!\!\big ]\big) \big )= \frac{ m d!}{\pi^d} \nu(\mathbb B_{d}).$$ Also, for the  ellipsoid $\mathbb B_{2,1}$, we have numerical evidence for such an equality. Taking all of this into account, we make the following conjecture. 
\begin{conj}
Suppose that the commuting $d$- tuple of operators $\boldsymbol{T}=(T_1,\ldots, T_d)$ is in the class $BS_{m,\vartheta}(\Omega)$. Then 
\[\text{\rm trace}\,\big (\text{dEt}\,\big(\big[\!\! \big [\boldsymbol{T}^*, \boldsymbol{T}\big ]\!\!\big ]\big) \big )\leq \frac{ m d!}{\pi^d} \nu(\overline{\Omega}), \]
where $\nu$ is the Lebesgue measure.
\end{conj}
We have made some progress towards finding a solution to this conjecture. 
\subsubsection*{Acknowledgment} The authors thank Dr. Cherian Varughese for several hours of fruitful discussions during the preparation of this manuscript.

\end{document}